\documentclass[mnsc,nonblindrev]{informs3_hide}
\OneAndAHalfSpacedXI 

\usepackage[english]{babel}
\usepackage[autostyle, english = american]{csquotes}
\MakeOuterQuote{"}
\usepackage[colorlinks,linkcolor=blue, citecolor=blue]{hyperref}

\usepackage{bm}
\usepackage{cleveref}
\crefname{subsection}{section}{subsections}

\newtheorem{observation}{Observation}

\usepackage{subfigure}

\newcommand{\tO}{\tilde{O}}

\usepackage{algorithm,algpseudocode}

\usepackage[dvipsnames]{xcolor}

\newcommand{\blue}{\textcolor{blue}}

\usepackage{xparse}

\NewDocumentEnvironment{myproof}{o}
{\IfNoValueTF{#1}{\paragraph{{Proof.} }} {\paragraph{{#1.} }} }
{\hfill$\Halmos$}

\usepackage{natbib,bbm,multirow,multicol}
 \bibpunct[, ]{(}{)}{,}{a}{}{,}%
 \def\bibfont{\small}%
\TheoremsNumberedThrough     
\ECRepeatTheorems

\EquationsNumberedThrough    

\begin{document}

\RUNAUTHOR{Chen, Jiang, Zhang, Zhou}
\RUNTITLE{Learning for Lost Sales Inventory Systems}

\TITLE{Learning to Order for Inventory Systems with Lost Sales and Uncertain Supplies}

\ARTICLEAUTHORS{
\AUTHOR{Boxiao Chen}
\AFF{College of Business Administration, University of Illinois at Chicago, Chicago, IL 60607,
\EMAIL{bbchen@uic.edu}}
\vspace{.05in}
\AUTHOR{Jiashuo Jiang}
\AFF{Stern School of Business, New York University, New York, NY 10012,
\EMAIL{jj2398@stern.nyu.edu}}
\vspace{.05in}
\AUTHOR{Jiawei Zhang}
\AFF{Stern School of Business, New York University, New York, NY 10012,
\EMAIL{jz31@stern.nyu.edu}}
\vspace{.05in}
\AUTHOR{Zhengyuan Zhou}
\AFF{Stern School of Business, New York University, New York, NY 10012,
\EMAIL{zzhou@stern.nyu.edu}}
}

\ABSTRACT{We consider a stochastic lost-sales inventory control system with lead time $L$ over a planning horizon $T$. Supply is uncertain, and is a function of the order quantity (due to random yield/capacity, etc). We aim to minimize the $T$-period cost, a problem that is known to be computationally intractable even under known distributions of demand and supply. In this paper, we assume that both the demand and supply distributions are unknown and develop a computationally efficient online learning algorithm. We show that our algorithm achieves a regret (i.e. the performance gap between the cost of our algorithm and that of an optimal policy over $T$ periods) of $\tO(L+\sqrt{T})$ when $L\geq\Omega(\log T)$. We do so by 1) showing our algorithm’s cost is higher by at most $\tO(L+\sqrt{T})$ for any $L\geq 0$ compared to an optimal constant-order policy under complete information (a well-known and widely-used algorithm) and 2) leveraging its known performance guarantee from the existing literature. To the best of our knowledge, a finite-sample $\tO(\sqrt{T})$ (and polynomial in $L$) regret bound when benchmarked against an optimal policy is not known before in the online inventory control literature.

A key challenge in this learning problem is that both demand and supply data can be censored; hence only truncated values are observable. We circumvent this challenge by showing that the data generated under an order quantity $q^2$ allows us to simulate the performance of not only $q^2$ but also $q^1$ for all $q^1<q^2$, a key observation to obtain sufficient information even under data censoring. By establishing a high probability coupling argument, we are able to evaluate and compare the performance of different order policies at their steady state within a finite time horizon. Since the problem lacks convexity, commonly used learning algorithms such as SGD and bisection cannot be applied, and instead, we develop an active elimination method that adaptively rules out suboptimal solutions.
}
\KEYWORDS{lost sales, lead time, supply uncertainty, online learning, censored data}

\maketitle
\section{Introduction}\label{sec:introduction}
Matching supply with demand is one of the key concepts in supply chain management. However, achieving this is not easy because of uncertainties in the system. One source of uncertainties stems from the randomness on the demand side. Consider an inventory control system, if demand realizes to be higher than the inventory level, there will be left-over inventories, and if demand realizes to be lower than the inventory level, some customer needs cannot be fulfilled. Based on the behavior of unsatisfied customers, typically the inventory system is classified as either backlogging or lost-sales \citep{devalve2020primal, bernstein2016simple, lei2020joint}. The lost-sales system is usually more relevant than the backlogging system, especially in the retailing setting where customers can easily go to a competitor when facing a stockout \citep{bijvank2011lost}. However, the lost-sales system is in general much harder to analyze than the backlogging counterpart, due to more complicated system dynamics, loss of convexity, etc. One classic model is the lost-sales inventory system with positive lead times, that is, it takes multiple periods for an order to arrive once it is placed. It is well-known that the optimal solution for this problem is computationally intractable, because the dimension of the state space of the underlying Markov decision process (MDP) equals the length of the lead time, which is also known as the curse of dimensionality \citep{bu2020constant}. Therefore, instead of striving to solve for the optimal solution, researchers shift their focus to developing effective heuristic policies \citep{levi20082, huh2009asymptotic, xin2016optimality, bu2020constant}.

Another source of uncertainties resides in the supply side as a result of the probabilistic nature of machine/capacity availability, yield, quality, and processing times \citep{yano1995lot, angkiriwang2014managing}. During and post the COVID-19 pandemic, supply uncertainties became more prominent worldwide, and supply chain runners need to adjust their operational strategies to better respond to supply disruptions, material shortages, and long lead times \citep{raj2022supply}.
The two most studied models on supply uncertainties are the proportional random yield model and the random capacity model. The former specifies that only a random proportion of the ordering quantity is fulfilled \citep{henig1990structure}, and the latter states that the ordering quantity is capped by a random capacity value \citep{ciarallo1994periodic}. Models with random yield are in general very hard to optimize, therefore, several heuristics were introduced in the literature and their performances were analyzed. See \citet{bu2020constant} for more discussions on the challenges of the random yield and random capacity problems.

Uncertainties from both the demand and the supply side contribute significantly to the difficulties of managing supply chains. In this paper, we consider both uncertainties. We study a periodic inventory control model over a planning horizon of $T$ periods. At the beginning of every period $t$, the company determines an ordering quantity $q_t$, of which only a portion $s(q_t, Z_t)$ will be fulfilled. The lead time is $L$, i.e., $s(q_t, Z_t)$ will be delivered after $L$ periods. Here, $s(\cdot, \cdot)$ is referred to as the supply function, and $Z_t$, $t=1, \ldots, T$, are independent and identically distributed (i.i.d.) random variables. As illustrated in Section \ref{sec:randomsupply},  $s(q_t, Z_t)$ considered in this paper include the random yield function $s(q_t, Z_t)=q_t Z_t$, the random capacity function $s(q_t, Z_t)=\min\{q_t, Z_t\}$, and so on. Note that when $s(q_t, Z_t)=q_t$, this is the classic lost-sales inventory system with lead times $L$ and deterministic supply, which is a special case of our models. Demand for period $t$, $D_t$, $t=1, \ldots, T$, are i.i.d. random variables. Demand $D_t$ is satisfied as much as possible by available inventories during period $t$. If demand is smaller than supply, left-over inventories will be carried over to the beginning of period $t+1$, for which a per-unit holding cost $h$ is charged. If demand realizes to be higher than supply, unfulfilled demand will be lost and a per-unit shortage penalty cost of $b$ will be incurred. Both the distributions of $D_t$ and $Z_t$ are unknown, and the company needs to infer their information only from historical data. The company would like to learn the demand and supply distributions and at the same time minimize the $T$ period total cost.

Even under complete information, that is, both the distributions of $D_t$ and $Z_t$ are known to the company, solving for the optimal solution is incredibly hard due to the presence of lost sales, lead times, and supply uncertainty. As noted earlier in the discussion, even when $s(q_t, Z_t)=q_t$, this problem is already computationally intractable. The only {asymptotically optimal} heuristic policy reported in the literature is the constant-order policy, which prescribes a constant quantity to be ordered every period independent of the starting state. \citet{bu2020constant} proves that the performance gap between the optimal constant-order policy and the true optimal policy (1) decays exponentially fast in the lead time $L$, (2) converges to $0$ when the penalty cost $b$ is large. In this paper, we consider a situation where neither the demand nor the supply distribution is known, and we develop online learning algorithms to learn the optimal constant-order policy using historical data.

Our main contributions are summarized as below.

\begin{enumerate}
\item Approaching the optimal policy with provable convergence rates. {In particular, we consider two benchmarks, one is the general optimal policy and the other is the optimal constant-order policy with given distribution information, which has been known to be asymptotic optimal. For any $L\geq0$ and any $T$,} we prove that the cost incurred by our learning algorithm is higher than that of the optimal constant-order policy by at most $\tO(L+\sqrt{T})$. On the other hand, it is shown in \citet{bu2020constant} that the cost of the optimal constant-order policy converges to the optimal policy exponentially fast in the lead time $L$, implying that the cost of our learning algorithm is higher than the \emph{optimal policy} by at most $\tO(L+\sqrt{T})$ when $L \ge \Omega(\log T)$. This is the first learning algorithm that provably approaches the optimal policy. In fact, even for the special case of lost-sales inventory system with lead times and deterministic supply, which has been extensively studied in the learning literature, there is no existing learning algorithm that approaches the optimal policy under any parameter regime. As will be discussed in Section \ref{sec:literature}, algorithms developed for this special case in \citet{huh2009adaptive}, \citet{zhang2020closing}, \citet{agrawal2022learning} and \citet{lyu2021ucb} are designed to approach various heuristic policies, and the best regret convergence rate in terms of its dependence in $L$ and $T$ is $\tO(L\sqrt{T})$, also benchmarked against certain heuristic policy. Our regret rate of $\tO(L+\sqrt{T})$ dominates the rate of $\tO(L\sqrt{T})$, not to mention that our regret also holds with supply uncertainty, and when $L \ge \Omega(\log T)$, our regret is benchmarked against the optimal policy.

\item Learning demand and supply distributions using censored data. We learn distributions for both the demand $D_t$ and the supply $Z_t$, departing from the existing learning literature that assumes the supply is deterministic. Both demand and supply data can be censored. Demand data is truncated by inventory levels, and the company can only observe the sales data instead of the true demand data. When supply has a random capacity, i.e., $s(q_t, Z_t)=\min\{q_t, Z_t\}$, the company can only observe the realized supply $s(q_t, Z_t)$ instead of the capacity value $Z_t$, i.e., the capacity data is truncated by the ordering decision. Because data is censored, it is not efficient to estimate the distributions of $D_t$ and $Z_t$ directly. We circumvent the data censoring issue using the following observation: for any two constant-order quantities $q^1<q^2$, utilizing the censored supply and demand data generated under $q^2$, we can construct approximate pseudo-costs to evaluate the performance of not only $q^2$ but also $q^1$ (Observation \ref{observe}, Lemma \ref{lem:computeIa}). Based on this critical observation, we develop simulation-based algorithms that simulate the performance of all order quantities after implementing the largest reasonable order quantity. This approach enables us to reduce the time spent on exploration and quickly focus on near optimal solutions.

\item \label{contribution3} Maximizing in an inventory system without convexity. In the literature of inventory control with online learning, a common property researchers rely on is that the objective function or part of it is convex in the decision variables. Based on convexity, some popular approaches can be applied, such as stochastic gradient decent (SGD) used in \citet{huh2009nonparametric, huh2009adaptive, zhang2018perishable, zhang2020closing, chen2019tailored}, and bisection used in \citet{agrawal2022learning, chen2019tailored, chen2020optimal}. However, convexity does not hold in our problem due to the complex structure of the random supply functions, therefore, commonly adopted approaches cannot be applied in our setting. Instead, we develop an adaptive elimination procedure that keeps eliminating suboptimal ordering values using historical censored data and maintains a shrinking active set. By letting the procedure to proceed in \textit{exponentially increasing} time intervals, we show that values that remain in the active set will perform very close to the true optimal constant-order policy with high probability as more data accumulates.

\item Applying high probability coupling arguments to evaluate policy performance at steady state. In order to solve for the optimal constant-order policy, we need to evaluate and compare the performance of different ordering policies at their steady state. However, this task can be highly nontrivial, because it may take a long time for a MDP to converge to its steady state and we only have a finite number of periods. Using a stochastic coupling argument, we prove that the MDP reaches its steady state after only $O(\log T)$ periods (Lemma \ref{lem:Costcoupling}) with high probability. Based on this result, we are able to adequately explore the inventory space and approach the optimal constant-order policy at a fast speed.
\end{enumerate}

\section{Literature Review} \label{sec:literature}
Supply chain management with supply uncertainties has been widely studied in the literature. See \citet{yano1995lot} and \citet{feng2018supply} for a detailed review of the area. As discussed earlier, one well-studied model is the proportional random yield model \citep{henig1990structure, kazaz2004production, federgruen2008selecting, federgruen2009optimal, li2013supply}. Models with random yield are in general very hard to solve for the optimal solution, therefore, papers such as \citet{bollapragada1999myopic, huh2010linear, inderfurth2015exact} develop heuristics and demonstrate their performance. A special case of the proportional random yield model is the all-or-nothing supply model \citep{anupindi1993diversification, tomlin2006value, babich2007competition, yang2009supply, gumucs2012supply}, based on which the random proportion takes either $0$ or $1$. Another commonly used model is the random capacity model. The random capacity model is first introduced in \citet{ciarallo1994periodic}, then extended to more general settings in papers such as \citet{wang1996periodic, chao2008joint, feng2010integrating} and \citet{bu2020constant}. All of the above-mentioned papers assume the distribution of the random supply is known, and none of them consider learning.

There is a stream of literature that studies the lost-sales inventory system with positive lead times and no supply uncertainties, that is $s(q_t, Z_t)=q_t$. Even for this simpler problem, as discussed earlier in the paper, the optimal solution cannot be directly solved for. Well-known heuristics developed for this problem include the base-stock policy \citep{janakiraman2004lost, huh2009asymptotic}, the capped base-stock policy \citep{xin2021understanding}, and the constant-order policy. \citet{reiman2004new} and \citet{zipkin2008old} first put forth the constant-order policy for the continuous review and periodic review lost-sales inventory system with positive lead times, respectively, and show that it performs favorably compared with other heuristics. \citet{goldberg2016asymptotic} applies the constant-order policy to the finite horizon lost-sales model and proves that it is asymptoticly optimal with large lead times. \citet{xin2016optimality} studies the constant-order policy for the infinite horizon lost-sales model and proves that the performance gap between it and the true optimal policy decays exponentially fast in the lead time. The constant-order policy is then generalized to the joint pricing and inventory control problem with lead times by \citet{chen2019asymptotic} and to several MDP settings in \citet{bai2020asymptotic}.

For the lost-sales inventory system with lead times and random supply functions that is considered in this paper, the constant-order policy is the only reported heuristic in the literature, and \citet{bu2020constant} proves that its performance converges to that of the true optimal policy exponentially fast in the lead time and is asymptotically optimal with large lost-sales penalty cost. All the above results regarding the constant-order policy assume the demand distribution is known, and there exist no existing algorithms that learn the constant-order policy when the demand distribution is unknown.

The area of inventory control with online demand learning has been flourishing in recent years \citep{lim2006model}. However, all existing studies assume supply is deterministic and demand is the only source of uncertainties. Algorithms are then proposed to learn only the demand distribution, which cannot be directly applied to the case when supply also has uncertainties. A special case of our problem is the well-known lost-sales inventory system with positive lead times $L$, but with deterministic supply. Because this special case is already too complex to have tractable optimal solutions, researchers propose online learning algorithms to learn various heuristics. \citet{huh2009adaptive} propose a gradient based learning algorithm that converges to the optimal base-stock heuristic policy. Their results are then improved by the SGD based learning algorithm in \citet{zhang2020closing} whose cost is proved to be higher than the optimal base-stock heuristic policy by at most $O(\exp(L) \sqrt{T})$. \citet{agrawal2022learning} develops a bisection based algorithm and further improves the result to $O(L\sqrt{T})$. By developing a UCB based learning algorithm for discrete demand, \cite{lyu2021ucb} proves that the cost of their algorithm is higher than the optimal capped base-stock heuristic policy by at most $O(L\sqrt{T})$. Different from the existing results, when $L \ge \Omega(\log T)$, the regret of our algorithm, $\tO(L+\sqrt{T})$, is obtained by comparing to the \emph{true optimal policy}. Other inventory control models with online demand learning include \cite{lin2022data} considering the repeated newsvendor problem and characterizing how the regret depends on a separation assumption about the demand, \citet{huh2009nonparametric} considering the lost-sales inventory system with zero lead time, \citet{zhang2018perishable} considering perishable products, \cite{chen2019tailored} studying the dual sourcing inventory system, \citet{yuan2021marrying} exploring inventory control problems with fixed setup cost. These works all develop SGD based learning algorithms and achieve a regret rate of $O(\sqrt{T})$. Problems with joint pricing and inventory decisions are explored in \citet{chen2019coordinating} and \citet{chen2021nonparametric} with regret $O(\sqrt{T})$.

{Finally, there has been a notable connection between reinforcement learning and inventory management. While inventory management under an online learning environment can be formulated as a reinforcement learning problem, a notable feature is that both the number of actions (the order quantity) and the number of states (the inventory in the pipeline) are actually infinite. In contrast, the algorithms developed for general reinforcement learning problem would require both the number of actions and the number of states (or aggregate states) to be finite and the regret bound would depend the numbers (e.g. \cite{jin2018q, dong2019provably}). Therefore, it is important to develop specialized algorithm for the inventory system separately (e.g. \cite{gong2021bandits} \cite{agrawal2022learning}). Our work falls into this line and compared with previous literature, we are the first to achieve the $\tO(L+\sqrt{T})$ regret convergence rate in the presence of positive lead time $L$.}

\section{Problem Formulation}
Consider a periodic-review lost-sales inventory system of a single product over a finite horizon of $T$ periods. The demand at each period $t$ is denoted by $D_t$, which belongs to the interval $[0,\bar{D}]$ and is assumed to be drawn independently from an unknown distribution $F(\cdot)$. At each period $t$, the company places an order of $q_t$, which will arrive after $L$ (a positive integer) periods. We consider the case where the company may not receive exactly what it orders. To model the supply uncertainty, we introduce a supply function $s(q,z): \mathbb{R}^2\rightarrow\mathbb{R}$, and we assume that the company at period $t$ receives a quantity given by $s(q_{t-L}, Z_t)$, where $Z_1,\dots, Z_T$ are i.i.d. non-negative random variables with a common distribution function $G(\cdot)$, which is assumed to belong to the interval $[\underline{\alpha}, \bar{\alpha}]$ and is assumed to be unknown. We also denote by $h$ the per-unit holding cost and denote by $b$ the per-unit lost-sale penalty cost, {which are assumed to be given at the very beginning}. We have the following sequence of events happening at each period $t$:\\
1. At the beginning of period $t$, the company observes the on-hand inventory level denoted by $I_t$ and all the inventories in pipeline ordered from the supplier, denoted by $(x_{1,t}, x_{2,t}, \dots, x_{L,t})$ where $x_{i,t}$ is the order quantity placed at period $t-L+i-1$ for $i=1,\dots,L$. The system state is $(I_t, x_{1,t}, x_{2,t}, \dots, x_{L,t})$.\\
2. The inventory order placed $L$ periods ago arrives and the random variable $Z_t$ is realized. Then, the on-hand inventory is increased to $I_t+s(x_{1,t}, Z_t)$.\\
3. The company placed an order with amount $q_t$ that will arrive at the beginning of period $t+L$.\\
4. The demand $D_t$ is realized and is satisfied as much as possible by the on-hand inventory. The unsatisfied demand is assumed to be lost and unobservable.\\
The objective of the company is to minimize the cumulative holding and penalty costs. The system state is updated as follows:
\[
I_{t+1}=(I_t+s(x_{1,t},Z_t)-D_t)^+,~~ x_{i,t+1}=x_{i+1,t}~~\forall 1\leq i\leq L-1\text{~~and~~}x_{L,t+1}=q_t
\]
A policy $\pi$ for the company is specified by the order quantities $q_1^{\pi},\dots, q_T^{\pi}$. Since the lost demand is assumed to be unobserved, we assume that only the \textit{censored demand} is known by the company, i.e., the company can only observe the sales quantity $\min\{I_t+s(x_{1,t},Z_t), D_t\}$ instead of the realization of $D_t$, and when $I_{t+1}=0$, the company does not know the volume of lost sales. Note that the \textit{supply data can be censored} as well, since only $s(x_{1,t}, Z_t)$ can be observed rather than $Z_t$, and $s(x_{1,t}, Z_t)$ may only contain truncated information about $Z_t$ (see more discussions on this issue in Section \ref{sec:randomsupply}). A policy $\pi$ is \textit{feasible} if and only if $\pi$ is non-anticipative, i.e, for each $t$, $q^{\pi}_t$ can only depend on the system state $(I_{\tau}^{\pi}, x^{\pi}_{1,\tau}, \dots, x^{\pi}_{L,\tau})$ for $\tau\leq t$ and the realized values of supply $(s_1,\dots,s_t)$. Note that the distribution functions $F(\cdot)$ and $G(\cdot)$ are assumed to be unknown by the company and need to be learned on-the-fly. Then, the cost incurred at period $t$ for the policy $\pi$ is denoted by
\[
C^{\pi}_t=h\cdot(I_t+s(x_{1,t},Z_t)-D_t)^++b\cdot (D_t-I_t-s(x_{1,t},Z_t))^+
\]
and the expected cumulative cost for the policy $\pi$ is denoted by
\begin{equation}\label{eqn:cumulativecost}
C^{\pi}(T, L)=\sum_{t=1}^{T}\mathbb{E}[C^{\pi}_t]=\sum_{t=1}^{T}\mathbb{E}\left[h\cdot(I_t+s(x_{1,t},Z_t)-D_t)^++b\cdot (D_t-I_t-s(x_{1,t},Z_t))^+\right]
\end{equation}
where $T$ is used to indicate the dependency on the number of periods in the entire horizon and $L$ is used to indicate the dependency on the lead time. Following this notation, the long-term average cost of the policy $\pi$ is denoted by
\begin{equation}\label{eqn:averagecost}
C^{\pi}_{\infty}=\limsup_{T\rightarrow\infty}\frac{1}{T}\cdot C^{\pi}(T,L)
\end{equation}
Following the standard conditions \citep{bu2020constant}, we will assume that the initial inventory is $I_1=0$ and the initial pipeline is also $0$, i.e., $x_{i,1}=0$ for all $1\leq i\leq L$.

\subsection{Constant Order Policies and Notion of Regret}
The optimal policy for minimizing the long-term reward in \eqref{eqn:averagecost} is known to be very complex and computationally intractable due to the curse of dimensionality caused by the lead time $L$. Thus, heuristics have been developed to solve the problem approximately. In this section, we introduce the heuristics studied in this paper, namely the \textit{constant order policies}, where the company places the same order in every period, regardless of the system state.

When the demand distribution and the supply function are unknown to the company, the optimal order quantity $q^*$ for minimizing \eqref{eqn:averagecost} cannot be directly computed. Our goal is to develop a feasible learning algorithm $\pi$. Using the optimal constant order policy $\pi_{q^*}$ as the benchmark, we measure the performance of the learning algorithm $\pi$ using the following notion of regret:
\begin{equation}\label{def:regret}
\text{Regret}^\pi_T=C^{\pi}(T, L)-T\cdot C^{\pi_{q^*}}_{\infty}
\end{equation}
An alternative way to define regret of online policy $\pi$ is to measure the additive difference between $C^{\pi}(T, L)$ and $C^{\pi_{q^*}}(T,L)$. We remark that for each policy $\pi$, the alternative regret will be at the same order of the regret defined in \eqref{def:regret} by noting that the gap between $C^{\pi_{q^*}}(T,L)$ and $T\cdot C^{\pi_{q^*}}_{\infty}$ can be bounded by $O(\sqrt{T})$ following standard concentration inequality for Markov chain with stationary distributions (see \Cref{lem:Markovconcentration}).

\subsection{Random Supply Function}\label{sec:randomsupply}
In this paper, we consider the random supply function $s(q,Z)$ that takes one of the following four formulations, {which is exogenously given to the company}:
\begin{enumerate}
    \item $s(q,Z)=q\cdot Z$.\label{form:1}
    \item $s(q,Z)=\min\{q,Z\}$.\label{form:2}
    \item $s(q,Z)=qZ/(q+\alpha Z^\rho)$ for $\rho\leq 1$ and $\alpha>0$.\label{form:3}
    \item $s(q,Z)=(qk)/(q+Z)$, for some $k>0$. \label{form:4}
\end{enumerate}
Formulation \ref{form:1} and \ref{form:2} covers the well-known random yield model and the random capacity model. Formulation \ref{form:3} is introduced in \citet{dada2007newsvendor} to model a non-linear relationship between the order quantity and the supply, which covers an increasing concave relation
of the output to the input over a wide range of parameters. Formulation \ref{form:4} has been used in \cite{cachon2003supply, tang2014pay} to study a situation where the supplier serves multiple firms and allocates the total output quantity, denoted by $k$, proportional to the firms’ order quantities, denoted by $q$. Note that the firm is not able to observe the order quantities required by other firms, which is captured by the random variable $Z$.

The above four formulations have been studied in \citet{feng2018supply}, which proves that all these four formulations are \textit{stochastically linear in mid-point} (Definition 1 in \citet{feng2018supply}). It has been shown in \cite{bu2020constant} that the long-run average cost of the optimal policy converges to the long-run average cost of the optimal constant order policy as $L\rightarrow\infty$, with the gap decreasing exponentially in the lead time $L$. This result justifies the efficiency of the constant order policy when the lead time $L$ is large, and is also the reason why we set the optimal constant order policy as the benchmark in the definition of regret in \eqref{def:regret}, which we further discuss in \Cref{remark:constant}.

Note that in formulations \ref{form:1}, \ref{form:3}, and \ref{form:4}, after observing $s(q,Z)$, the value of $Z$ can be inferred. However, this does not hold for formulation \ref{form:2}, where the value of $Z$ is truncated by the ordering quantity $q$. That is, if $Z$ realizes to be higher than $q$, then the company can only observe $q$. This data censoring issue for supply uncertainty creates extra challenges for estimating the distribution of $Z$.

The following observation plays a critical role in addressing the supply data censoring issue, and it is a key step to develop our learning algorithm (further explained in \Cref{sec:simulation}).
\begin{observation}\label{observe}
If the random supply function takes one of the Formulation \ref{form:1}, \ref{form:2}, \ref{form:3} and \ref{form:4}, then for any $q$ and $Z$, as long as we observe the value of $s(q,Z)$, we know the value of $s(q',Z)$ for any $q'\leq q$.
\end{observation}
Clearly, for formulations \ref{form:1}, \ref{form:3} and \ref{form:4}, this observation holds true by noting that the value of $Z$ can actually be derived backward from the value of $q$ and $s(q,Z)$. For formulation \ref{form:2}, $q=s(q,Z)$ implies $q\leq Z$. Then, for any $q'\leq q$, we must have $s(q',Z)=q'$. Also, $q>s(q,Z)$ implies $q>Z=s(q,Z)$. Then, for any $q'\leq q$, we have $s(q',Z)=\min\{q',s(q,Z)\}$. Therefore, we justify \Cref{observe} also holds for formulation \ref{form:2}.

\begin{remark}
Note that when the supply is deterministic, i.e., $s(q, Z)=q$, it is a special case of our models. This is the classic lost-sales inventory system with positive lead times that is extensively studied in the literature \citep{bijvank2011lost}.
\end{remark}

\begin{remark}
{Note that the above four formulations satisfy the \textit{stochastic linearity in mid-point} condition proposed in \cite{feng2018supply} and it follows that $C^{\pi_q}_{\infty}$ can be transferred into a convex function of $\mu=\mathbb{E}_{Z\sim G}[s(q,Z)]$. Therefore, bounding regret \eqref{def:regret} admits a convex reformulation if we regard $\mu$ as the decision variable. However, in our problem, the actual decision is $q$. Then, it would require additional efforts to learn the correspondence relationship between $q$ and $\mu$ since distribution $G$ is unknown and only censored feedback can be observed. Instead, as shown in the next section, our approach does not require convexity and the only property we need is \Cref{observe}. Therefore, we can deal with other random supply function formulations where \Cref{observe} holds but stochastic linearity with mid-point property fails to hold, and our main result \Cref{thm:mainregret} still applies.}
\end{remark}

\section{Algorithm and General Description of Our Approach}
In this section, we propose our learning algorithm to achieve the regret of optimal order. We begin with re-formulating the long-run average cost of a constant order policy $\pi_q$. Note that in expression \eqref{eqn:cumulativecost}, the true value of $D_t$ is unobservable due to lost sales and censored demand. However, we now show that in order to learn the optimal order quantity $q^*$, it is enough to focus solely on the on-hand inventory $I_t$ and the supply $s(x_{1,t}, Z_t)$, which can be directly observed.

First, under the constant order policy $\pi_q$, the on-hand inventory is updated as follows:
\begin{equation}\label{eqn:ConstantQ}
  I^{\pi_q}_{t+1}=(I_t^{\pi_q}+s(q,Z_t)-D_t)^+.
\end{equation}
From queueing theory \citet{asmussen2008applied}, the sequence $\{I^{\pi_q}_t\}_{t=1}^{\infty}$ converges in probability to a random variable $I^{\pi_q}_{\infty}$, which we refer to as the limiting inventory level under the constant order policy $\pi_q$, as long as the following condition is satisfied for the order quantity $q$:
\begin{equation}\label{eqn:conditiononQ}
  \mathbb{E}_{Z\sim G}[s(q,Z)]< \mathbb{E}_{D\sim F}[D].
\end{equation}
If condition in (\ref{eqn:conditiononQ}) is not met, {the on-hand inventory level will approach infinity in the long run}. Therefore, we also impose the same condition for our analyses as shown in Assumption \ref{assump:OrderUB} below.
\begin{assumption}\label{assump:OrderUB}
The company knows an upper bound of the optimal order quantity $q^*$, denoted by $\bar{q}$, that satisfies $\mathbb{E}[s(\bar{q},Z)]<\mathbb{E}[D]$.
\end{assumption}
Assumption \ref{assump:OrderUB} is very mild, because any value not satisfying the condition in (\ref{eqn:conditiononQ}) will cause { the on-hand inventory level to approach infinity}, therefore those values can be easily detected as suboptimal.

We have
\[
C^{\pi_q}_{\infty}=h\cdot\mathbb{E}\left[(I^{\pi_q}_{\infty}+s(q,Z)-D)^+\right] +b\cdot\mathbb{E}\left[(D-s(q,Z)-I^{\pi_q}_{\infty})^+\right].
\]
Moreover, it holds that
\[
I^{\pi_q}_{\infty}=^d (I^{\pi_q}_{\infty}+s(q,Z)-D  )^+
\]
where $=^d$ denotes identical in distribution. By taking expectation over both sides of the following equation,
\[
I^{\pi_q}_{\infty}+s(q,Z)-D=[I^{\pi_q}_{\infty}+s(q,Z)-D]^+-[D-s(q,Z)-I^{\pi_q}_{\infty}]^+,
\]
we have that
\[\begin{aligned}
\mathbb{E}[I^{\pi_q}_{\infty}]+\mathbb{E}[s(q,Z)]-\mathbb{E}[D]&=\mathbb{E}[I^{\pi_q}_{\infty}+s(q,Z)-D]^+-\mathbb{E}[D-s(q,Z)-I^{\pi_q}_{\infty}]^+\\
&=\mathbb{E}[I^{\pi_q}_{\infty}]-\mathbb{E}[D-s(q,Z)-I^{\pi_q}_{\infty}]^+
\end{aligned}\]
which implies that
\begin{equation}\label{eqn:averageConstantQ}
  C^{\pi_q}_{\infty}=h\cdot\mathbb{E}[I^{\pi_q}_{\infty}]+b\cdot \mathbb{E}[D]-b\cdot\mathbb{E}[s(q,Z)].
\end{equation}
Note that $C^{\pi_q}_{\infty}$ is unobservable, because the term $\mathbb{E}[D]$ in \eqref{eqn:averageConstantQ} is unobservable due to demand censoring. However, the term $\mathbb{E}[D]$ is independent of the order quantity $q$. In order to obtain the optimal order quantity $q^*$, it is equivalent to minimize the \textit{pseudo-cost} defined as follows:
\begin{equation}\label{eqn:psuedocost}
  \hat{C}^{\pi_q}_{\infty}=h\cdot\mathbb{E}[I^{\pi_q}_{\infty}]-b\cdot\mathbb{E}[s(q,Z)]
\end{equation}
over the set $Q=\{q: q \le \bar q\}$. Here, the pseudo-cost $\hat{C}^{\pi_q}_{\infty}$ is observable. However, $\hat{C}^{\pi_q}_{\infty}$ is not convex in the order quantity $q$, in which case commonly used learning approaches such as SGD and bisection cannot be applied. For more discussions regarding this issue, see Section \ref{sec:introduction} point \ref{contribution3}.

We now describe our learning algorithm for solving \eqref{eqn:psuedocost}. Speaking at a high level, we transfer our problem into a multi-arm bandit problem by specifying $K+1$ points uniformly over the interval $[0,\bar{q}]$ i.e., we specify a set of points $\mathcal{A}=\{a_1,\dots,a_{K+1}\}$ such that $a_k=\frac{k-1}{K}\cdot\bar{q}$ for any $k=1,\dots,K+1$. Then, our algorithm proceeds in epochs $n=1,2,\dots$ by maintaining an active set $\mathcal{A}_n\subset\mathcal{A}$ for each epoch $n$. The key element of our algorithm is to guarantee that for each epoch $n$ and each point $a\in\mathcal{A}_n$, the gap between $\hat{C}^{\pi_a}_{\infty}$ and $\hat{C}^{\pi_{q^*}}_{\infty}$ is upper bounded by $\gamma_n$, where $\{\gamma_n\}_{n\geq1}$ is a decreasing sequence to be determined later. To be specific, we let each epoch $n$ contain $\max\{\frac{1}{\gamma_{n+1}^2}\cdot\log T, 3L\}$ number of time periods and the implementation of our algorithm at epoch $n$ can be classified into the following three steps:\\
1. We implement the constant order policy $\pi_{a^{n*}}$, where $a^{n*}$ is the largest element in the active set $\mathcal{A}_n$. \\
2. We use the censored demand to \textit{simulate} the pseudo-cost of the policies $\pi_a$ for each $a\in\mathcal{A}_n$ and we construct a confidence interval of $\hat{C}^{\pi_a}_{\infty}$ for each $a\in\mathcal{A}_n$ (simulation step further discussed in \Cref{sec:simulation}).\\
3. We use the constructed confidence intervals to identify $\mathcal{A}_{n+1}\subset\mathcal{A}_n$ such that for each element $a\in\mathcal{A}_{n+1}$, the gap between $\hat{C}^{\pi_a}_{\infty}$ and $\hat{C}^{\pi_{a^*}}_{\infty}$ is upper bounded by $(h+b)\cdot\gamma_{n}$, where $\hat{C}^{\pi_{a^*}}_{\infty}=\min_{a\in\mathcal{A}}\hat{C}^{\pi_{a}}_{\infty}$.

Following the steps outlined above, as our learning algorithm proceeds and $n$ increases, the active set $\mathcal{A}_n$ shrinks and the optimal order quantity $q^*$ is gradually approximated. Our algorithm is formally described in \Cref{alg}. Note that the implementation of \Cref{alg} depends on a fixed constant $\kappa_2$. We provide further discussion on how to select $\kappa_2$ in \Cref{sec:diskappa}. By specifying the value of $K$ and the sequence $\{\gamma_n\}_{n\geq1}$, we are able to prove the following theorem regarding the regret upper bound of our algorithm, which is the main theorem of our paper.
\begin{algorithm}
\caption{Learning-based Constant Order Policy}\label{alg}
\begin{algorithmic}[1]
\State \textbf{Input:} $K$ and $\{\gamma_n\}_{n\geq1}$.
\State Initialize $\mathcal{A}_1=\mathcal{A}$, where $\mathcal{A}=\{a_1,\dots,a_{K+1}\}$ such that $a_k=\frac{k-1}{K}\cdot\bar{q}$ for any $k=1,\dots,K+1$.
\State Set $\tau_n=\sum_{n'=1}^{n-1}\kappa_2\cdot\max\{\frac{1}{\gamma_{n'+1}^2}\cdot\log T, 3L\}+1$ as the start of epoch $n$ for each $n\geq1$, where $\kappa_2$ is a fixed constant.
\For {epoch $n=1,2,\dots,$}
\State Identify $a^{n*}$ as the largest element in the active set $\mathcal{A}_n$.
\For {time period $t=\tau_n$ to $\tau_{n+1}-1$}
\State Implement the constant order policy $\pi_{a^{n*}}$.
\State Observe the value of the supply $s(x_{1,t}, Z_t)$ and the on-hand inventory level $I_t$.
\EndFor
\State For each $a\in\mathcal{A}_n$, we construct $\tilde{C}^a_n$ as follows:
\begin{itemize}
  \item obtain the simulated supply $s(a, Z_t)$ under policy $\pi_a$ for each $t=\tau_n+L,\dots,\tau_{n+1}-1$;
  \item starting from $I_{\tau_n+L}^a=I_{\tau_n+L}$, for $t=\tau_n+L,\dots,\tau_{n+1}-1$, do the following:
   \begin{itemize}
     \item if $I_{t+1}>0$, then $I^a_{t+1}=(I^a_t+s(a,Z_t)+I_{t+1}-I_t-s(a^{n*},Z_t))^+$;
     \item if $I_{t+1}\leq0$, then $I^a_{t+1}=0$.
   \end{itemize}
  \item compute
  \[\begin{aligned}
  \tilde{C}^a_n=&h\cdot\frac{1}{\tau_{n+1}-\tau_n-\kappa_2\max\{\log T,2L\}}\cdot\sum_{t=\tau_n+\kappa_2\max\{\log T,2L\}}^{\tau_{n+1}-1}I_t^a\\
  &-b\cdot\frac{1}{\tau_{n+1}-\tau_n-\kappa_2\max\{\log T,2L\}}\cdot\sum_{t=\tau_n+\kappa_2\max\{\log T,2L\}}^{\tau_{n+1}-1}s(a,Z_t).
  \end{aligned}    \]
\end{itemize}
\label{step:simulate}
\State Denote by $\tilde{C}^*_n=\min_{a\in\mathcal{A}_n} \tilde{C}^a_n$ and identify the active set for epoch $n+1$.
\begin{equation}\label{eqn:activesetevolve}
\mathcal{A}_{n+1}=\{ a\in\mathcal{A}_n: \tilde{C}^a_n\leq \tilde{C}^*_n+(h+b)\cdot\frac{\gamma_n}{2} \}
\end{equation}
\EndFor
\end{algorithmic}
\end{algorithm}
\begin{theorem}\label{thm:mainregret}
Denote by $\pi$ \Cref{alg} with input $K=\sqrt{T}$ and $\gamma_n=2^{-n}$ for each $n\geq1$. Suppose that the random supply function takes one of the four formulations specified in \Cref{sec:randomsupply}. Then, under \Cref{assump:OrderUB}, the regret of $\pi$ has the following upper bound:
\[
\text{Regret}(\pi)\leq \kappa\cdot\kappa_2\cdot(L+\sqrt{T})\cdot\log T
\]
where $\kappa$ is a constant that is independent of $L$ and $T$, and $\kappa_2$ is the constant used in \Cref{alg}.
\end{theorem}
\begin{remark}\label{remark:constant}
We note that \Cref{thm:mainregret} implies a regret bound of \Cref{alg} even compared to the \textit{optimal policy}, when the lead time $L$ is sufficiently large. To see this, we apply Theorem 1 in \citet{bu2020constant} to show that $C^{\pi_{q^*}}_{\infty}-C^{\pi^*}_{\infty}\leq \kappa_3\cdot\gamma^{L}$, where $\kappa_3$ and $\gamma\in(0,1)$ are constants and $\pi^*$ stands for the optimal policy, i.e. $\pi^*=\text{argmin}_{\pi} C^{\pi}_{\infty}$. Therefore, we have $C^{\pi}(T, L)- T\cdot C^{\pi^*}_{\infty}\leq \kappa\cdot\kappa_2\cdot(L+\sqrt{T})\cdot\log T+\kappa_3\cdot T\cdot\gamma^{L}$, which implies that $C^{\pi}(T, L)- T\cdot C^{\pi^*}_{\infty}\leq \tO(L+\sqrt{T})$ when $L\geq \Omega(\log T)$. This is the \textit{first} time that a sublinear regret bound is derived for an online policy with respect to the optimal policy.  
As a result, our result justifies the efficiency of constant order policies for inventory control systems with large lead time, under an online learning environment.
\end{remark}

In the literature, the most related results on regret convergence rates are derived from \citet{huh2009adaptive}, \citet{zhang2020closing}, \citet{agrawal2022learning} and \citet{lyu2021ucb}, which study a special case of our problem with deterministic supply. The state-of-the-art regret convergence rate is $O(L\sqrt{T})$, derived in \citet{agrawal2022learning} for continuous demand benchmarked against the optimal base-stock heuristic policy and \citet{lyu2021ucb} for discrete demand benchmarked against the optimal capped base-stock heuristic policy. Our regret rate of $\tO(L+\sqrt{T})$ compares favorably with the existing results in this special case in terms of the dependence on $L$ and $T$, and it is derived benchmarked against the optimal policy (instead of a heuristic policy) when $L \ge \Omega(\log T)$.

{We now discuss the tightness of the $\tO(L+\sqrt{T})$ bound shown in \Cref{thm:mainregret}. On the one hand, from Proposition 1 of \cite{zhang2020closing}, we know that no learning algorithm can achieve a regret bound better than $\Omega(\sqrt{T})$ when compared to the optimal policy. On the other hand, since the lead time is $L$ and any adjustment over the order quantity will only influence the system after $L$ periods, we know that no learning algorithm can achieve a regret better than $\Omega(L)$. Therefore, when compared to the optimal policy, no learning algorithm can achieve a regret better than $\Omega(L+\sqrt{T})$, which implies our regret bound is tight up to a logarithmic term when $L\geq\Omega(\log T)$. }

\begin{remark}
{Regarding the improvement from $\tilde{O}(L\sqrt{T})$ to $\tilde{O}(L+\sqrt{T})$, indeed both the choice of benchmarks and the algorithmic design/analysis matter. Please note one key difference over the confidence intervals established in our paper and in \cite{agrawal2022learning} and \cite{lyu2021ucb} is that our confidence interval is \textit{independent} of the lead time $L$, while the confidence intervals in \cite{agrawal2022learning} and \cite{lyu2021ucb} depends on $L$. This property is made possible by the specific structure of the constant-order policy benchmark. Then, by further exploiting this structure, we can show that during each epoch, {we only incur a loss $O(L+\sqrt{N_i})$ with $N_i$ being the number of periods in epoch $i$}, while the loss in \cite{agrawal2022learning} and \cite{lyu2021ucb} scales similar to $O(L\cdot\sqrt{N_i})$. Therefore, by using an exponentially increasing epoch, we guarantee the number of epochs is at most $O(\log T)$, {which leads to our final regret bound $\tilde{O}(L+\sqrt{T})$}.}
\end{remark}

\begin{remark}
\Cref{alg} can be carried out efficiently. Note that for each quantity $a\in\mathcal{A}$, the inventory level $I^a_t$ is simulated at most once for each period $t=1,\dots,T$ and $|\mathcal{A}|=\sqrt{T}$. Therefore, the overall computation complexity of \Cref{alg} is upper bounded by $O(T^{\frac{3}{2}})$.
\end{remark}

\subsection{Discussion on the Simulation Step}\label{sec:simulation}
In this section, we discuss why we could use the censored demand of the constant order policy $\pi_{a^{n*}}$ to simulate the pseudo-cost of the policies $\pi_a$ for each $a\in\mathcal{A}_n$, as outlined in step 10 in \Cref{alg}. Following \Cref{observe}, since $a^{n*}$ is the largest element in the active set $\mathcal{A}_n$, after observing the value of $s(a^{n*}, Z_t)$, we know the value of $s(a, Z_t)$ for all $a\in\mathcal{A}_n$, for any $t=\tau_n+L,\dots,\tau_{n+1}-1$. Thus, we can use $s(a, Z_t)$ for any $t=\tau_n+L,\dots,\tau_{n+1}-1$ to approximate the term $\mathbb{E}[s(a,Z)]$ in the expression \eqref{eqn:psuedocost} for $\hat{C}^{\pi_a}_{\infty}$, for all $a\in\mathcal{A}_n$. Following Hoeffding's inequality, the approximation error can be bounded with a high probability (formalized in \Cref{sec:confidenceinterval}).

For any $a\in\mathcal{A}_n$, we can approximate $\mathbb{E}[I^{\pi_a}_{\infty}]$. We define a stochastic process $\{I^a_t\}_{t=\tau_n+L}^{\tau_{n+1}-1}$ revolving in the following way:
\begin{equation}\label{eqn:aevolve}
  I^a_{\tau_n+L}=I_{\tau_n+L}\text{~and~}I^a_{t+1}=(I^a_t+s(a,Z_t)-D_t)^+\text{~for~all~}t=\tau_n+L,\dots,\tau_{n+1}-2
\end{equation}
Clearly, when the value of $D_t$ is censored, we can not directly obtain the value of $I^a_{t+1}$. However, we now show that if the random supply function takes one of the four formulations specified in \Cref{sec:randomsupply}, we can use the on-hand inventory level $I_t$ to derive the value of $I^a_{t}$, for any $t=\tau_n+L+1,\dots,\tau_{n+1}-1$. Note that $\{I_t\}_{t=\tau_n+L}^{\tau_{n+1}-1}$ evolves in the following way:
\begin{equation}\label{eqn:Itevolve}
  I_{t+1}=(I_t-s(a^{n*},Z_t)-D_t)^+.
\end{equation}
Suppose that the value of $I^a_t$ is known, we derive the value of $I^a_{t+1}$ under the following two cases.\\
1. If $I_{t+1}>0$, then from \eqref{eqn:Itevolve}, we can obtain the value of $D_t$ and we can derive the value of $I^a_{t+1}$ directly following \eqref{eqn:aevolve}.\\
2. If $I_{t+1}\leq0$, then we have $I_t\leq D_t-s(a^{n*},Z_t)$. Note that $s(a,Z_{\tau})\leq s(a^{n*}, Z_{\tau})$ for all $\tau\leq t$, we must have $I^a_t\leq I_t$. Thus, we have
\[
I^a_t\leq I_t\leq D_t-s(a^{n*},Z_t)\leq D_t-s(a,Z_t)
\]
which implies that $I^a_{t+1}=0$.

The above two steps are formalized in the following lemma.
\begin{lemma}\label{lem:computeIa}
Suppose that the stochastic process $\{I^a_t\}_{t=\tau_n+L}^{\tau_{n+1}-1}$ is defined in \eqref{eqn:aevolve} and denote by $\{I_t\}_{t=\tau_n+L}^{\tau_{n+1}-1}$ the on-hand inventory level evolving in \eqref{eqn:Itevolve}. Then, the value of $I^a_t$ can be computed iteratively for $t=\tau_n+L, \dots,\tau_{n+1}-2$ in the following way:
\begin{itemize}
  \item if $I_{t+1}>0$, then $I^a_{t+1}=(I^a_t+s(a,Z_t)+I_{t+1}-I_t-s(a^{n*},Z_t))^+$;
  \item if $I_{t+1}\leq0$, then $I^a_{t+1}=0$.
\end{itemize}
\end{lemma}
After deriving the value of $\{I^a_t\}_{t=\tau_n+L}^{\tau_{n+1}-1}$, we use this sequence to approximate $\mathbb{E}[I^{\pi_a}_{\infty}]$. The key is to establish the coupling between the stochastic process $\{I^a_t\}_{t=\tau_n+L}^{\tau_{n+1}-1}$ and another stochastic process, which we further explain in \Cref{sec:confidenceinterval}.

{We do note that for \blue{supply function} 1, 3, 4 specified in \Cref{sec:randomsupply}, we can directly obtain the value of $Z_t$ and we know $s(q',Z_t)$ for any $q'$. However, knowing $s(q',Z_t)$ itself is not enough to simulate the value $I^{q'}_t$ since the demand is \textit{censored} as well. In contrast, as discussed above, choosing the largest element in $\mathcal{A}_n$ would allow us to simulate $I^a_t$ for all other $a\in\mathcal{A}_n$. In fact, choosing the largest element is one key \blue{distinguishing feature} between our algorithm and the classical active arm elimination algorithm in the multi arm bandit (MAB) literature \citep{even2006action}. In MAB literature, each element in the active set will be chosen a certain number of times. However, if we directly apply this algorithm in our problem, we need to change the order quantity $|\mathcal{A}_n|$ times in each epoch $n$. \blue{Note that each time the order quantity is changed}, it incurs a $O(L)$ loss because of the lead time. Therefore this would lead to a worse dependency on $L$, such as $O(L\cdot \sqrt{N_i})$ where $N_i$ denotes the number of periods in each epoch. The key innovation of our algorithm is that we can learn everything we need (censored supply and censored demand) by simply choosing the largest element of the active set in each epoch. Since the total number of epochs is $O(\log T)$, as shown later in \Cref{sec:pfThm}, there will only be a loss $O(L\cdot\log T)$ incurred from changing the order quantity.}

\subsection{Discussion on the constant $\kappa_2$}\label{sec:diskappa}
Note that the implementation of \Cref{alg} depends on a fixed constant $\kappa_2$. We now discuss how should we select the value of $\kappa_2$.

In order for the regret bound in \Cref{thm:mainregret} to hold, a condition on the constant $\kappa_2$ would be $\kappa_2\geq \delta(F,G,\bar{q})$, where $\delta(F,G,\bar{q})$ is a constant that depends solely on $F, G$ and $\bar{q}$, and is independent of $L$ and $T$. Though the value of $\delta(F,G,\bar{q})$ is unknown at the beginning since we assume the distributions $F$ and $G$ are unknown, we can simply set $\kappa_2=\log T$ and the condition $\kappa_2\geq \delta(F,G,\bar{q})$ will automatically be satisfied when $T$ is large enough. Such an operation will only induce an additional multiplicative $\log T$ term into the final regret bound in \Cref{thm:mainregret}. Another way is to spend the first $O(\sqrt{T})$ periods as a pure learning phase to learn the distributions $F$ and $G$, and estimate an upper bound of $\delta(F,G,\bar{q})$, which is a constant independent of $T$ and $L$. Such an operation will only induce an additional additive $O(\sqrt{T})$ term into the final regret bound in \Cref{thm:mainregret}, which arises from the learning phase.

\section{Proof of Regret Bound}\label{sec:pfThm}
In this section, we prove the regret bound in \Cref{thm:mainregret}. Our analysis can be classified into the following four steps:\\
1) we establish the Lipschitz continuity of the pseudo-cost $\hat{C}^{\pi_q}_{\infty}$ over $q$. As a result, instead of comparing with $\hat{C}^{\pi_{q^*}}_{\infty}$, we can compare with $\hat{C}^{\pi_{a^*}}_{\infty}$ where $a^*=\text{argmin}_{a\in\mathcal{A}} \hat{C}^{\pi_{a}}_{\infty}$. We show that the additional regret term caused by this replacement of benchmark is at most $O(\sqrt{T})$.\\
2) we provide a bound over the gap between the actual pseudo-cost incurred at each epoch $n$ and the long-term average $\hat{C}^{\pi_{a^{n*}}}_{\infty}$. The proof of the bound relies on a novel coupling argument between two stochastic processes, which is explained in \Cref{sec:gap}.\\
3) we denote by $\mathcal{E}$ the event that for each epoch $n$ (except the last epoch), the pseudo-cost of each $a\in\mathcal{A}_n$ falls into the confidence interval $[\tilde{C}^a_n-(h+b)\cdot\frac{\gamma_n}{2}, \tilde{C}^a_n+(h+b)\cdot\frac{\gamma_n}{2}]$, i.e.,
\begin{equation}\label{eqn:eventE}
  \mathcal{E}=\{ |\tilde{C}^a_n-\hat{C}^{\pi_a}_{\infty}|\leq(h+b)\cdot\frac{\gamma_n}{2},~\forall a\in \mathcal{A}_n, \forall 1\leq n\leq N-1 \}
\end{equation}
where $N$ denotes the total number of epochs. We show that event $\mathcal{E}$ occurs with a high probability.\\
4) we show how $a^*$ is approximated by the revolution of the active set $\mathcal{A}_n$ in \eqref{eqn:activesetevolve}, which leads to our final regret bound.

Following the above four steps, we decompose the regret of our policy $\pi$ as follows:
\begin{align}
\text{Regret}(\pi)=&\sum_{n}\sum_{t=\tau_n}^{\tau_{n+1}}(h\cdot\mathbb{E}[I^{\pi}_t]-b\cdot\mathbb{E}[s(a^{n*},Z_t)]-\hat{C}_{\infty}^{\pi_{q^*}}) \label{eqn:regretdecompose}\\
=&\underbrace{\sum_{t=1}^T(\hat{C}^{\pi_{a^*}}_{\infty}-\hat{C}^{\pi_{q^*}}_{\infty})}_{I}+ \underbrace{\sum_{n}\sum_{t=\tau_n}^{\tau_{n+1}}(h\cdot\mathbb{E}[I^{\pi}_t]-b\cdot\mathbb{E}[s(a^{n*},Z_t)]-\hat{C}_{\infty}^{\pi_{a^{n*}}})}_{II}   \nonumber\\
&+\underbrace{\sum_{n}\sum_{t=\tau_n}^{\tau_{n+1}} (\hat{C}^{\pi_{a^{n*}}}_{\infty}-\hat{C}^{\pi_{a^*}}_{\infty})}_{III} \nonumber
\end{align}
We use the Lipschitz continuity established in the first step to bound the term I in \eqref{eqn:regretdecompose}. We use the high probability bound established in the second step to bound the term II in \eqref{eqn:regretdecompose}. We finally use step three and step four to bound the term III in \eqref{eqn:regretdecompose}.
\subsection{Proof of Lipschitz Continuity}\label{sec:Lipschitz}
In this section, we establish the Lipschitz continuity of $\mathbb{E}[I^{\pi_q}_{\infty}]$ over order quantity $q$. We denote by $\hat{s}(\mu,Z)=s(q,Z)$ for $q$ satisfying $\mathbb{E}[s(q,Z)]=\mu$. Our approach relies on existing result (\Cref{lem:convex}) showing that if we interpret the psuedo-cost as a function over $\mu$, then this function is a convex function, which implies Lipschitz continuity since $\mu$ belongs to a bounded region. Moreover, for the random supply function taking one of the four formulations specified in \Cref{sec:randomsupply}, one can check that if we interpret $\mu$ as a function of $q$, then this function is Lipschitz continuous. Therefore, we prove the Lipschitz continuity of $\mathbb{E}[I^{\pi_q}_{\infty}]$. We summarize our result in the following lemma, where the proof is relegated to \Cref{sec:pfap}.
\begin{lemma}\label{lem:Lipschitz}
There exists a constant $\beta>0$ such that for any $q_1, q_2\in[0,\bar{q}]$, we have
\[
|\hat{C}^{\pi_{q_1}}_{\infty}-\hat{C}^{\pi_{q_2}}_{\infty}|\leq\beta\cdot|q_1-q_2|.
\]
\end{lemma}

\subsection{Gap Between Actual Pseudo Cost and Long-term Average Pseudo Cost}\label{sec:gap}
We provide the bound over the gap between the actual pseudo cost incurred during each epoch $n$ and the pseudo-cost $\hat{C}^{\pi_{a^{n*}}}_{\infty}$. Our proof relies on establishing the coupling between the stochastic process $\{I_t\}_{t=\tau_n}^{\tau_{n+1}-1}$ and the stochastic process defined as follows:
\begin{equation}\label{eqn:tIevolve}
  \tilde{I}^{a^{n*}}_{\tau_n}=^d I^{\pi_{a^{n*}}}_{\infty}\text{~and~}\tilde{I}^{a^{n*}}_{t+1}=(\tilde{I}^{a^{n*}}_t+s(a^{n*},Z_t)-D_t)^+\text{~for~}t=\tau_n,\dots,\tau_{n+1}-2
\end{equation}
It is clear to see that the distribution of $\tilde{I}^{a^{n*}}_t$ is identical to the distribution of $I^{\pi_{a^{n*}}}_{\infty}$, for each $t=\tau_n,\dots,\tau_{n+1}-1$. The coupling argument is formalized in the following lemma.
\begin{lemma}\label{lem:Costcoupling}
Denote by $N$ the total number of epochs and denote by $\mathcal{B}$ the event that $I_{\tau_n}\leq \kappa_1\cdot\log T$ and $\tilde{I}^{a^{n*}}_{\tau_n}\leq \kappa_1\cdot\log T$, where $\kappa_1>0$ is a fixed constant, and $\{I_{\tau_n+\kappa_2\cdot\max\{\log T, 2L\}}=\tilde{I}^{a^{n*}}_{\tau_n+\kappa_2\cdot\max\{\log T, 2L\}}\}$, for every epoch $n\in[N]$, i.e.,
\[
\mathcal{B}=\{ I_{\tau_n}\leq \kappa_1\cdot\log T, \tilde{I}_{\tau_n}^{a^{n*}}\leq\kappa_1\cdot\log T\text{~and~}\{I_{\tau_n+\kappa_2\cdot\max\{\log T, 2L\}}=\tilde{I}^{a^{n*}}_{\tau_n+\kappa_2\cdot\max\{\log T, 2L\}}\},~\forall n \}.
\]
Then, we have that
\[
P(\mathcal{B})\geq 1-\frac{3N}{T^2}.
\]
\end{lemma}
The proof is relegated to \Cref{sec:pfap}.

From \Cref{lem:Costcoupling}, we know that conditioning on the event $\mathcal{B}$, it holds that $I_t=\tilde{I}^{a^{n*}}_{t}$ for any $t=\tau_n+\kappa_2\cdot\max\{\log T, 2L\},\dots,\tau_{n+1}-1$ and any epoch $n\in[N]$. Moreover, note that the distribution of $\tilde{I}^{a^{n*}}_t$ is identical to the distribution of $I^{a^{n*}}_{\infty}$. It holds that
\begin{equation}\label{eqn:Reintercost}
\hat{C}^{\pi_{a^{n*}}}_{\infty}=h\cdot\mathbb{E}[\tilde{I}^{a^{n*}}_t]-b\cdot\mathbb{E}[s(a^{n*},Z_t)],~~\forall t=\tau_n,\dots,\tau_{n+1}-1,~\forall n\in[N]
\end{equation}
As a result, the expected value of $I_t$ for $t=\tau_n+\kappa_2\cdot\max\{\log T, 2L\},\dots,\tau_{n+1}-1$ will be the same as the expected value of $I^{a^{n*}}_{\infty}$, which implies that the expected actual cost should be the same as the long-term average cost for $t=\tau_n+\kappa_2\cdot\max\{\log T, 2L\},\dots,\tau_{n+1}-1$. Thus, we can obtain an upper bound over the gap between the actual pseudo cost and the long-term average pseudo cost $\hat{C}^{\pi_{a^{n*}}}_{\infty}$, for each epoch $n\in[N]$. By summing up the bound for each epoch $n\in[N]$, we get an upper bound of the term II in \eqref{eqn:regretdecompose} for the entire horizon, which is formalized in the following lemma.
\begin{lemma}\label{lem:Gapcost}
It holds that
\[
\sum_{n}\sum_{t=\tau_n}^{\tau_{n+1}}(h\cdot\mathbb{E}[I^{\pi}_t]-b\cdot\mathbb{E}[s(a^{n*},Z_t)]-\hat{C}_{\infty}^{\pi_{a^{n*}}})\leq hN\cdot \kappa_1\kappa_2\log T\cdot\max\{\log T,2L\}+3hN\bar{D}
\]
\end{lemma}
{Note that $\hat{C}_{\infty}^{\pi_{a^{n*}}}$ admits the formulation \eqref{eqn:Reintercost}. Therefore, it suffice to bound the gap between $\mathbb{E}[I^{\pi}_t]$ and $\mathbb{E}[\tilde{I}_t^{a^{n*}}]$.} We refer the detailed proof to \Cref{sec:pfap}.
\subsection{Probability Bound on the Event $\mathcal{E}$}\label{sec:confidenceinterval}
We now show that the pesudo-cost of each $a\in\mathcal{A}_n$ at each epoch $n$ falls into the confidence interval $[\tilde{C}^a_n-\gamma_n, \tilde{C}^a_n+\gamma_n]$ with a high probability and we provide a bound over the probability that event $\mathcal{E}$ happens. The key is to establish the stochastic coupling between the stochastic process $\{I^a_t\}_{t=\tau_n}^{\tau_{n+1}-1}$ defined in \eqref{eqn:aevolve} and the stochastic process $\{\tilde{I}^a_t\}_{t=\tau_n}^{\tau_{n+1}-1}$ defined as follows:
\begin{equation}\label{eqn:tildeIat}
  \tilde{I}^a_{\tau_n}=^d I^{\pi_a}_{\infty}\text{~and~}\tilde{I}^a_{t+1}=(\tilde{I}^a_t+s(a,Z_t)-D_t)^+ \text{~for~}t=\tau_n,\dots,\tau_{n+1}-2
\end{equation}
We formalize the coupling argument in the following lemma, which generalizes the stochastic coupling established in \Cref{lem:Costcoupling} from the implemented order quantity $a^{n*}$ to all quantity $a\in\mathcal{A}_n$.
\begin{lemma}\label{lem:confidencecoupling}
Denote by $N$ the total number of epochs and denote by $\mathcal{C}$ the event that $I^a_{\tau_n}\leq \kappa_1\cdot\log T$ and $\tilde{I}^a_{\tau_n}\leq\kappa_1\cdot\log T$, where $\kappa_1>0$ is a fixed constant, and $\{I^a_{\tau_n+\kappa_2\cdot\max\{\log T, 2L\}}=\tilde{I}^{a}_{\tau_n+\kappa_2\cdot\max\{\log T, 2L\}}\}$, for every epoch $n\in[N]$ and every $a\in\mathcal{A}_n$, i.e.,
\[
\mathcal{C}=\{ I^a_{\tau_n}\leq \kappa_1\cdot\log T, ~\tilde{I}_{\tau_n}^{a}\leq\kappa_1\cdot\log T\text{~and~}\{I^a_{\tau_n+\kappa_2\cdot\max\{\log T, 2L\}}=\tilde{I}^{a}_{\tau_n+\kappa_2\cdot\max\{\log T, 2L\}}\},~\forall n\in[N], \forall a\in\mathcal{A}_n \}
\]
Then, we have that
\[
P(\mathcal{C})\geq 1-\frac{3(K+1)N}{T^2}.
\]
where $K$ is given as the input of \Cref{alg} to denote $|\mathcal{A}|$.
\end{lemma}
The proof is relegated to \Cref{sec:pfap}.

For each epoch $n$ and each action $a\in\mathcal{A}_n$, it is clear to see that the distribution of $\tilde{I}_t^a$ is identical to the distribution of $I^{\pi_a}_{\infty}$. Therefore, we can use the average value of $\tilde{I}^a_t$ for $t=\tau_n+\kappa_2\cdot\max\{\log T,2L\}$ to $\tau_{n+1}-1$ to approximate the value of $\mathbb{E}[I^{\pi_a}_{\infty}]$, where the length of the confidence interval can be given by $\gamma_n$. Further note that conditioning on the event $\mathcal{C}$ happens, the value of $\{I^a_t\}_{t=\tau_n+\kappa_2\cdot\max\{\log T,2L\}}^{\tau_{n+1}-1}$ equals the value of $\{\tilde{I}^a_t\}_{t=\tau_n+\kappa_2\cdot\max\{\log T,2L\}}^{\tau_{n+1}-1}$, which implies that $\hat{C}^{\pi_a}_{\infty}\in[\tilde{C}^a_n-\gamma_n, \tilde{C}^a_n+\gamma_n]$ with a high probability.
\begin{lemma}\label{lem:highprob}
We have the following bound over the probability that event $\mathcal{E}$ happens, where event $\mathcal{E}$ is defined in \eqref{eqn:eventE},
\[
P(\mathcal{E})\geq 1-\frac{7(K+1)N}{T^2}.
\]
\end{lemma}
The proof is relegated to \Cref{sec:pfap}.

\subsection{Proof of \Cref{thm:mainregret}}
We are now ready to prove our main theorem. Following \eqref{eqn:regretdecompose}, we have
\begin{align}
\text{Regret}(\pi)=&\underbrace{\sum_{t=1}^T(\hat{C}^{\pi_{a^*}}_{\infty}-\hat{C}^{\pi_{q^*}}_{\infty})}_{I}+ \underbrace{\sum_{n}\sum_{t=\tau_n}^{\tau_{n+1}}(h\cdot\mathbb{E}[I^{\pi}_t]-b\cdot\mathbb{E}[s(a^{n*},Z_t)]-\hat{C}_{\infty}^{\pi_{a^{n*}}})}_{II}   \nonumber\\
&+\underbrace{\sum_{n}\sum_{t=\tau_n}^{\tau_{n+1}} (\hat{C}^{\pi_{a^{n*}}}_{\infty}-\hat{C}^{\pi_{a^*}}_{\infty})}_{III} \nonumber
\end{align}
We use the Lipschitz continuity established in \Cref{lem:Lipschitz} to bound the term I. We denote by $a'\in\mathcal{A}$ the nearest one to $q^*$. Clearly, from the construction of the set $\mathcal{A}$, we know that $|q^*-a'|\leq\frac{\bar{q}}{2K}$. Therefore, from \Cref{lem:Lipschitz}, we know that
\begin{equation}\label{eqn:boundI}
\text{I}=T\cdot (\hat{C}^{\pi_{a^*}}_{\infty}-\hat{C}^{\pi_{q^*}}_{\infty})\leq T\cdot (\hat{C}^{\pi_{a'}}_{\infty}-\hat{C}^{\pi_{q^*}}_{\infty})\leq T\cdot\frac{\beta\bar{q}}{2K}=\frac{\beta\bar{q}\sqrt{T}}{2}.
\end{equation}
where we note $K=\sqrt{T}$.

We now bound the term II.
From \Cref{lem:Gapcost}, we know that
\begin{equation}\label{eqn:boundII}
\text{II}=\sum_{n}\sum_{t=\tau_n}^{\tau_{n+1}}(h\cdot\mathbb{E}[I^{\pi}_t]-b\cdot\mathbb{E}[s(a^{n*},Z_t)]-\hat{C}_{\infty}^{\pi_{a^{n*}}})\leq hN\cdot \kappa_1\kappa_2\log T\cdot\max\{\log T,2L\}+3hN\bar{D}
\end{equation}
We now proceed to bound the term III with the help of the probability bound established in \Cref{sec:confidenceinterval}.

We now assume that the event
\[
\mathcal{E}=\{ |\tilde{C}^a_n-\hat{C}^{\pi_a}_{\infty}|\leq(h+b)\cdot\frac{\gamma_n}{2},~\forall a\in \mathcal{A}_n, \forall 1\leq n\leq N-1 \}
\]
happens. For each epoch $n$ and each $a\in\mathcal{A}_{n+1}$, from \eqref{eqn:activesetevolve} and the conditions of event $\mathcal{E}$, we have
\[
\hat{C}_{\infty}^{\pi_a}-\hat{C}_{\infty}^{\pi_{a^*}}\leq \tilde{C}^a_n-\tilde{C}_n^{a^*}+(h+b)\cdot\gamma_n\leq\frac{3}{2}\cdot(h+b)\cdot\gamma_n.
\]
Note that $a^{(n+1)*}\in\mathcal{A}_{n+1}$, we have that
\[
\hat{C}_{\infty}^{\pi_{a^{(n+1)*}}}-\hat{C}_{\infty}^{\pi_{a^*}}\leq \frac{3}{2}\cdot(h+b)\cdot\gamma_n.
\]
which implies the following inequality conditional on the event $\mathcal{E}$ happens,
\[
\text{III}=\sum_{n}\sum_{t=\tau_n}^{\tau_{n+1}} (\hat{C}^{\pi_{a^{n*}}}_{\infty}-\hat{C}^{\pi_{a^*}}_{\infty})\leq\frac{3(h+b)}{2}\cdot\sum_{n=1}^N\sum_{t=\tau_n}^{\tau_{n+1}}\gamma_{n-1}
\]
Moreover, denote by $N$ the total number of epochs. We have
\[
\kappa_2\cdot\sum_{n=1}^{N-1}\frac{1}{\gamma_n^2}\cdot\log T\leq \sum_{n=1}^{N-1}\kappa_2\cdot\max\{\frac{1}{\gamma_n^2}\cdot\log T, 3L\}\leq T
\]
which implies that
\[
\sum_{n=1}^{N-1}\frac{1}{\gamma_n^2}\leq\frac{T}{\kappa_2\cdot\log T}
\]
By specifying $\gamma_n=2^{-n}$, we have that $N\leq\log_4 \frac{3T+\log T}{\kappa_2\cdot\log T}$. Therefore, conditional on the event $\mathcal{E}$ happens, we have that
\begin{align}
\sum_{n=1}^{N}\sum_{t=\tau_n}^{\tau_{n+1}-1}\hat{C}_{\infty}^{a^{n*}}-\hat{C}_{\infty}^{a^*} &\leq\frac{3(h+b)}{2}\cdot\sum_{n=1}^{N}\sum_{t=\tau_n}^{\tau_{n+1}-1} \gamma_{n-1}=3(h+b)\cdot\sum_{n=1}^{N} \gamma_n\cdot\max\{\frac{1}{\gamma_n^2}\cdot\log T, 3L\}  \label{eqn:5181}\\
&=3(h+b)\cdot\sum_{n=1}^{\lfloor \log_4 \frac{3L}{\log T} \rfloor}\gamma_n\cdot 3L+3(h+b)\cdot\sum_{n=\lfloor \log_4 \frac{3L}{\log T} \rfloor+1}^{N}\frac{\log T}{\gamma_n}\nonumber\\
&\leq 3(h+b)\cdot\sum_{n=1}^{\lfloor \log_4 \frac{3L}{\log T} \rfloor}\gamma_n\cdot 3L+ 3(h+b)\cdot\sum_{n=1}^{N}\frac{\log T}{\gamma_n}\nonumber\\
&\leq 9(h+b)L+3(h+b)(2^{N+1}-1)\log T\nonumber\\
&\leq 9(h+b)L+6(h+b)\cdot\sqrt{\frac{(3T+\log T)\log T}{\kappa_2}}\nonumber
\end{align}
If the event $\mathcal{E}$ does not happen, clearly, we have that
\[
\text{III}=\sum_{n}\sum_{t=\tau_n}^{\tau_{n+1}} (\hat{C}^{\pi_{a^{n*}}}_{\infty}-\hat{C}^{\pi_{a^*}}_{\infty})\leq T\cdot(h+b)\cdot \bar{D}
\]
where we note that $\hat{C}^{\pi_{a^{n*}}}_{\infty}\leq(h+b)\cdot\bar{D}$ for each $n$. Therefore, we have the following upper bound over the term III,
\begin{align}
\text{III}&=\mathbb{E}\left[\sum_{n}\sum_{t=\tau_n}^{\tau_{n+1}} (\hat{C}^{\pi_{a^{n*}}}_{\infty}-\hat{C}^{\pi_{a^*}}_{\infty})\mid\mathcal{E}\right]\cdot P(\mathcal{E})+\mathbb{E}\left[\sum_{n}\sum_{t=\tau_n}^{\tau_{n+1}} (\hat{C}^{\pi_{a^{n*}}}_{\infty}-\hat{C}^{\pi_{a^*}}_{\infty})\mid\mathcal{E}^c\right]\cdot(1-P(\mathcal{E}))  \label{eqn:boundIII}\\
&\leq \mathbb{E}\left[\sum_{n}\sum_{t=\tau_n}^{\tau_{n+1}} (\hat{C}^{\pi_{a^{n*}}}_{\infty}-\hat{C}^{\pi_{a^*}}_{\infty})\mid\mathcal{E}\right]+T\cdot(h+b)\cdot \bar{D}\cdot (1-P(\mathcal{E})) \nonumber\\
&\leq 9(h+b)L+6(h+b)\cdot\sqrt{\frac{(3T+\log T)\log T}{\kappa_2}}+\frac{7(K+1)N\bar{D}(h+b)}{T}\nonumber
\end{align}
where the last inequality follows from \eqref{eqn:5181} and the probability bound on the event $\mathcal{E}$ from \Cref{lem:highprob}. Combining \eqref{eqn:boundI}, \eqref{eqn:boundII}, and \eqref{eqn:boundIII}, we have that
\[
\text{Regret}(\pi)= \text{I}+\text{II}+\text{III}\leq \kappa\cdot\kappa_2\cdot(L+\sqrt{T})\cdot\log T
\]
where $\kappa$ is a constant that is independent of $L$ and $T$.
Therefore, our proof of our main result \Cref{thm:mainregret} is completed.

{
\section{Numerical Experiments}
In this section, we numerically investigate the performance of our \Cref{alg}. We implement \Cref{alg} on a periodic-review lost-sale inventory systems and illustrate how the performance of \Cref{alg} would vary over the model parameters, i.e., critical ratio, random supply variance, lead time, and demand distribution variance. We measure the performance of \Cref{alg} by the relative average regret defined as
\begin{equation}\label{eqn:030901}
\frac{C^{\pi}(T,L)-C^{\pi_{q^*}}(T, L)}{C^{\pi_{q^*}}(T, L)}\times 100\%
\end{equation}
where $\pi$ denotes \Cref{alg} and $C^{\pi_{q^*}}(T, L)$ denotes the total cost of the optimal constant order policy. Note that following a standard concentration inequality for a Markov chain with stationary distributions (see \Cref{lem:Markovconcentration}), $C^{\pi_{q^*}}(T, L)$ would be close to $T\cdot C_{\infty}^{\pi_{q^*}}$. Therefore, the relative average regret \eqref{eqn:030901} reflects how well the performance of \Cref{alg} is when compared to the optimal constant order policy.

\subsection{Random Capacity Model}\label{sec:RCapacity}

In this section, we consider the random capacity model where $s(q,Z)=\min\{q, Z\}$.
In \Cref{fig:01}, we study how the critical ratio of the model would influence the performance of \Cref{alg}. \blue{The model parameters of the experiment are given as follows. Note that the distribution of $Z$ and the demand distribution is unknown to \Cref{alg}, however, it is assumed to be given for the optimal constant order policy and the optimal policy.} We let $Z$ follow a uniform distribution over $[5, 15]$. We also set the demand distribution to be a normal distribution with mean $10$ and variance $4$, truncated on $0$. We set the lead time $L=10$ and fix the holding cost $h=5$. We vary the value of $b$ to study the performance of \Cref{alg} under different critical ratio, which is defined as $b/(h+b)$. To be more concrete, we set $b=28.33$, $20$ and $15$ such that the critical ratio would equal to $0.85$, $0.8$ and $0.75$. As we can see, for all critical ratios, the average regret decreases with $T$ and finally converge to within $10\%$ when $T=1000$. \blue{Moreover, the regret increases with critical ratio. Indeed, by noting \eqref{eqn:activesetevolve}, a larger $b$ results in a larger active set $\mathcal{A}_n$ for each epoch $n$. This implies that the elimination of the non-optimal elements in the active set in our algorithm becomes slower when $b$ is large, which leads to a larger regret. Therefore, we see the regret increases with the critical ratio.}

\begin{figure}[!h]
    \centering
    \subfigure[The performance of \Cref{alg} with different critical ratios.]{
    \includegraphics[width=0.4\textwidth]{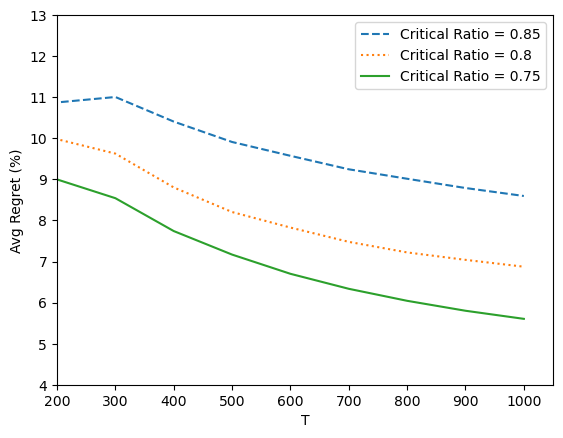}
    \label{fig:01}}
    \subfigure[The performance of \Cref{alg} with different supply variance.]{
    \includegraphics[width=0.4\textwidth]{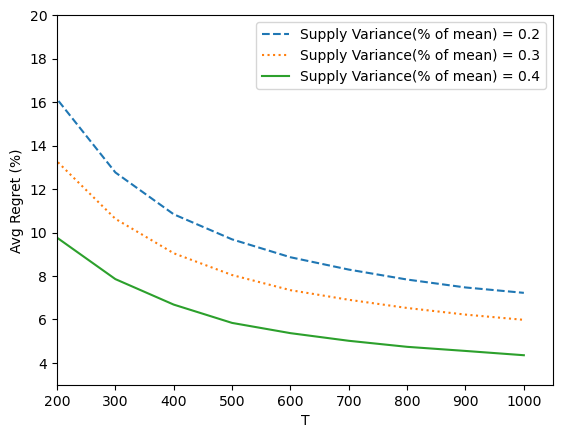}
    \label{fig:02}}
    \subfigure[The performance of \Cref{alg} with different lead time.]{
    \includegraphics[width=0.4\textwidth]{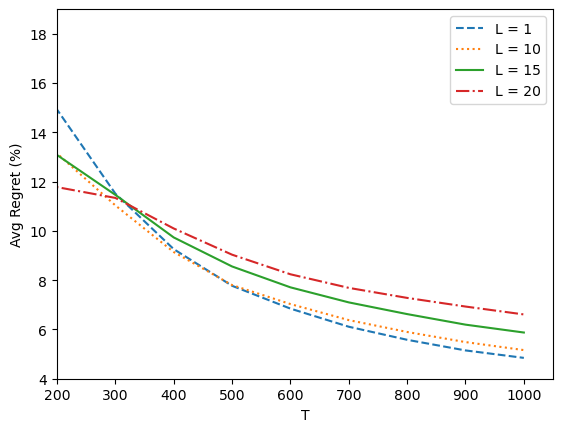}
    \label{fig:03}}
    \subfigure[The performance of \Cref{alg} with different demand variance.]{
    \includegraphics[width=0.4\textwidth]{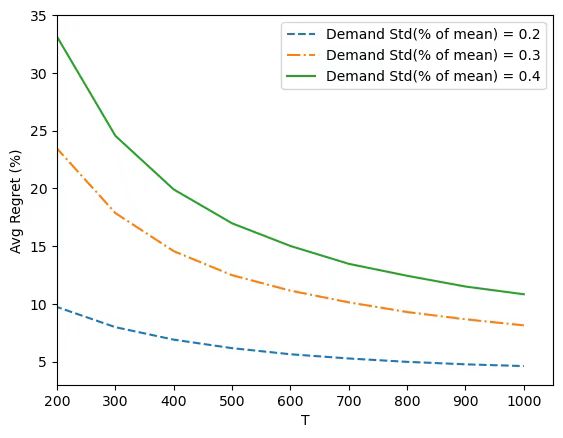}
    \label{fig:04}}
    \caption{Numerical results for different model parameters under random capacity supply.}
\end{figure}

In \Cref{fig:02}, we study how the variance of the supply parameter $Z$ would influence the performance of \Cref{alg}. We remain the same setting by letting $s(q,Z)=\min\{q, Z\}$, where $Z$ follows a uniform distribution over $[10-a, 10+a]$, with $a$ being set to $2$, $3$ and $4$. Again, we set the demand distribution to be a normal distribution with mean $10$ and variance $4$, truncated on $0$. We set the lead time $L=10$, the holding cost $h=5$ and the lost-sale penalty cost $b=5$. As we can see, \blue{the regret decreases with the supply variance. From our numerical observation, we find out one main source of regret comes from the first multiple epochs where we order a large quantity. However, when the supply variance increases, we observe that the optimal constant order quantity also increases, which actually make it easier for our algorithm to approximate the optimal constant order quantity since our algorithm always starts from the largest possible order quantity and gradually shrinks to the optimal one. Therefore, a larger supply variance results in a smaller regret incurred during the first multiple epochs, and as a result, a smaller regret during the entire horizon.}

In \Cref{fig:03}, we study how the performance of \Cref{alg} depends on the lead time. We consider the random capacity model by setting $s(q,Z)=\min\{q, Z\}$, where $Z$ follows a uniform distribution over $[5, 15]$. The demand follows a normal distribution with mean $10$ and variance $4$, truncated on $0$. We set the holding cost $h=5$ and the lost-sale penalty cost $b=5$. In general, the regret increases with the lead time $L$, which is in correspondence with our $\tO(L+\sqrt{T})$ regret. \blue{Note that from the numerical result, the regret decreases in the lead time when T is relatively small. This may due to when $T$ is relatively small, the number of epochs is small and especially for the first few epochs, the number of periods is linear in $L$ according to the design of our algorithm. However, when $L$ is quite small, e.g., $L=1$, the number of periods in the first few epochs is also very small and the estimation error of the expected cost is large. This leads to a large regret when $L$ is small. Therefore, when $T$ is relatively small, the reason why the regret decreases in the lead time relies on the fact that the estimation error is large when $L$ is small.} In \Cref{fig:04}, we set $L=10$ and study how the demand variance influences the performance of \Cref{alg}. As we can see, the regret {increases} as the demand variance becomes larger.

\subsection{Random Yield Model}\label{sec:RYield}
In this section, we consider the random yield model where $s(q,Z)=q\cdot Z$. The experiment setup is the same as the setup in \Cref{sec:RCapacity}. As we can see, the results for the random yield model are similar to the results for the random capacity model in \Cref{sec:RCapacity}. In \Cref{fig:05}, we study how the critical ratio of the model would influence the performance of \Cref{alg}. We fix the holding cost $h=5$. We set the value of $b$ to be $28.33$, $20$ and $15$ to study the performance of \Cref{alg} under different critical ratio. The regret \blue{increases} with critical ratio. In \Cref{fig:02}, we study how the variance of the supply parameter $Z$ influences the performance of \Cref{alg}. We let $Z$ follow a uniform distribution over $[10-a, 10+a]$, with $a$ being set to $2$, $3$ and $4$. As we can see, the influence of the supply variance on the regret is negligible and for all instances, the relative regret shrinkages to within $5\%$ as $T$ grows to $1000$. In \Cref{fig:03}, we show that the regret increases with the lead time $L$, which is in correspondence with our $\tO(L+\sqrt{T})$ regret. In \Cref{fig:04}, we show that the regret {increases} as the demand variance becomes larger.

\begin{figure}[!h]
    \centering
    \subfigure[The performance of \Cref{alg} with different critical ratios.]{
    \includegraphics[width=0.4\textwidth]{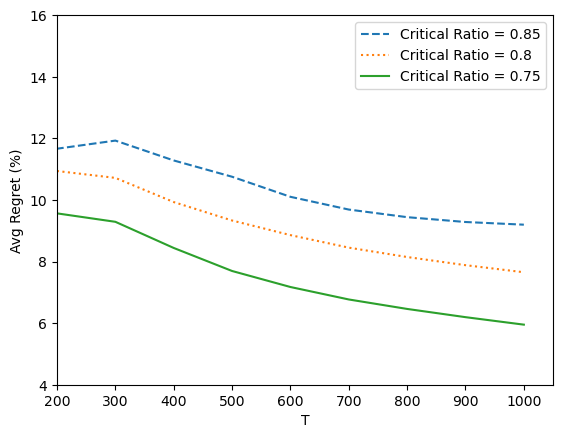}
    \label{fig:05}}
    \subfigure[The performance of \Cref{alg} with different supply variance.]{
    \includegraphics[width=0.4\textwidth]{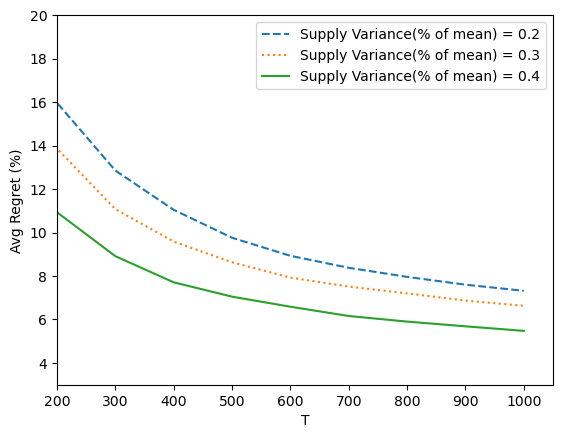}
    \label{fig:06}}
    \subfigure[The performance of \Cref{alg} with different lead time.]{
    \includegraphics[width=0.4\textwidth]{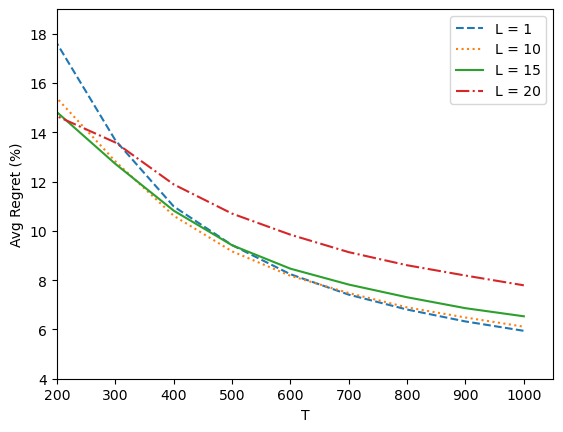}
    \label{fig:07}}
    \subfigure[The performance of \Cref{alg} with different demand variance.]{
    \includegraphics[width=0.4\textwidth]{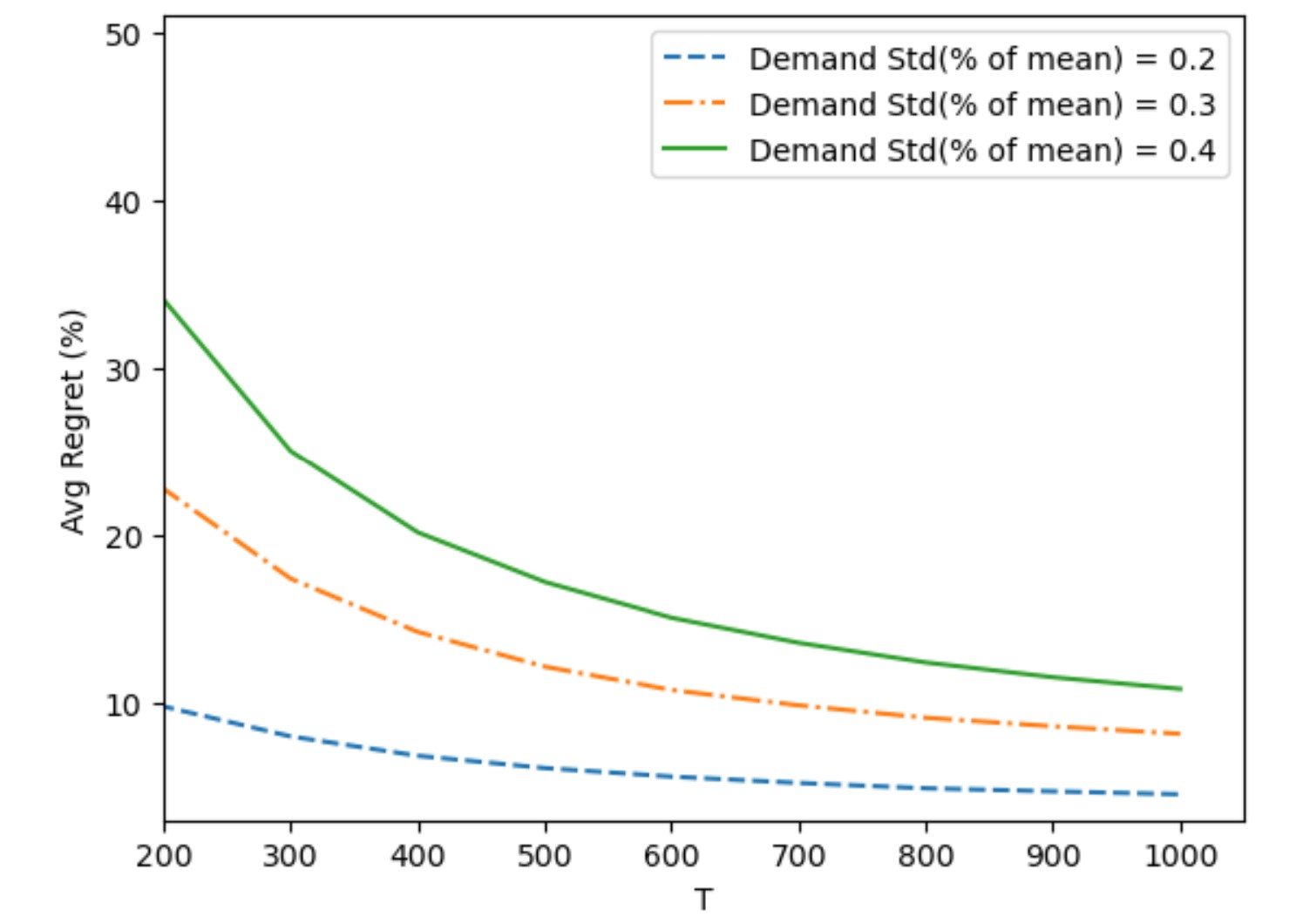}
    \label{fig:08}}
    \caption{Numerical results for different model parameters under random yield supply.}
\end{figure}

\subsection{Comparison with the Optimal Policy}
In this section, we compare the performance of \Cref{alg} to the optimal policy which can be computed by solving the dynamic program (DP). We also compare the performance of the optimal constant order policy with that of the DP to illustrate the effectiveness of the optimal constant order policy to serve as the benchmark. We use the random yield model with $T=1000$ and $b=h=5$. In order to compute the DP, we discretize the state space of $(I_t, x_{1,t}, x_{2,t}, \dots, x_{L,t})$ and the demand distribution is set to be a normal distribution with mean $10$ and variance $9$, discretized to the same magnitude as the state space. Due to the curse of dimensionality, the computation complexity for solving the DP would grow exponentially in the lead time $L$. Therefore, for computational tractability, we set the lead time to be $2$. The numerical results are reported in \Cref{fig:DP}. As we can see, the performance of the optimal constant order policy is quite close to that of the DP, with a relative difference within $2\%$. This difference is expected to become even smaller when the lead time $L$ increases, as illustrated by the theoretical findings of \cite{bu2020constant}.

\begin{figure}[!h]
  \centering
  \subfigure[The performance of optimal constant order policy with respect to the DP. The average regret is defined as $(\text{Optimal Constant Order}-\text{DP})/\text{DP}$]{
    \includegraphics[width=0.4\textwidth]{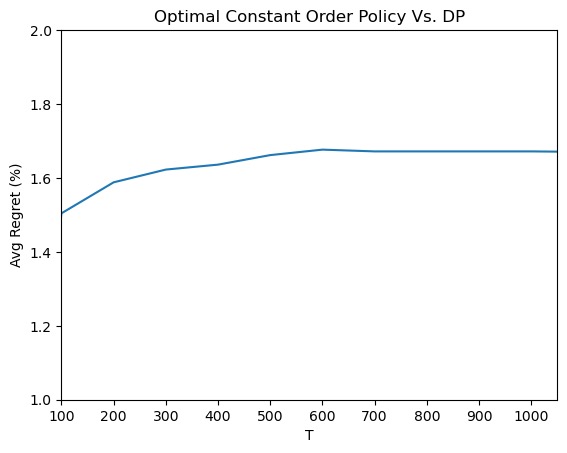}}
    \subfigure[The performance of \Cref{alg} with respect to the DP. The average regret is defined as $(\text{\Cref{alg}}-\text{DP})/\text{DP}$]{
    \includegraphics[width=0.4\textwidth]{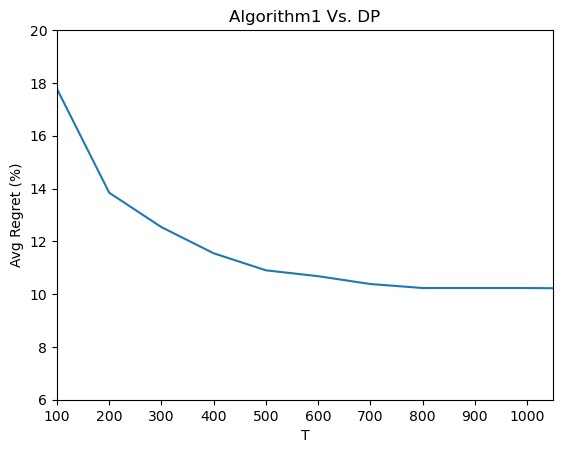}}
    \caption{Numerical results with respect to the optimal DP benchmark.}
    \label{fig:DP}
\end{figure}
}

\section{Conclusion}

In this paper, we study the lost-sales inventory system with lead times $L$. Both demand and supply have uncertainties, for which the distributions are unknown. This departs from the existing literature on online learning that assumes supply is deterministic and only considers demand uncertainty. The company needs to learn the demand and supply distributions from historical censored data. Demand censoring is caused by the fact that demand data is truncated by the inventory level, and supply censoring stems from the fact that capacity data is truncated by the ordering quantity. Because of demand and supply data censoring, it is not feasible to measure the performance of a policy directly, which requires the knowledge of the full demand and supply distributions. To circumvent this obstacle, we adopt a pseudo cost measure and prove that for any two constant-order quantities $q^1<q^2$, using the censored data generated under $q^2$, we can simulate the pseudo cost for not only $q^2$ but also $q^1$. The critical observation enables us to significantly reduce the time spent on exploration. In order to evaluate the performance of a policy under steady state, we develop a high probability coupling argument to show that the MDP under our policy approaches its steady state within $O(\log T)$ periods. Note that the objective function of our problem lacks convexity, a property that is utilized by almost all existing papers in the literature of inventory control with learning. We propose an active elimination based algorithm to achieve a regret of $\tO(L+\sqrt{T})$ when compared with the optimal constant-order policy, and when $L \ge \Omega(\log T)$, our algorithm approaches the optimal policy at the same rate.

There are many interesting directions for future research. For example, in the current setting, pricing is not considered. It would be a set of nice results if pricing can be included in the decision process and a learning algorithm can be developed accordingly. Another direction is to consider multiple products, where there exist substitution effects between different products and learning algorithms need to learn the substitution behavior of customers.

\bibliographystyle{apalike}
\bibliography{myreferences}
\clearpage

\begin{APPENDICES}
\section{Useful Previous Results}\label{sec:usefullemma}
We first state the well-known Hoeffdeing's inequality, which establishes concentration bound for i.i.d. random variables.
\begin{lemma}[Hoeffding's Inequality]\label{lem:Hoeffding}
Let $X_1,\dots,X_m$ be independent random variables such that $a_i\leq X_i\leq b_i$ almost surely for each $i\in[m]$. Then, denote by $S_n=\sum_{i=1}^m X_i$. It holds that
\[
P(S_n-\mathbb{E}[S_n]\geq t)\leq\exp(-\frac{2t^2}{\sum_{i=1}^m(b_i-a_i)^2})
\]
\end{lemma}
We then state a general concentration bound for Markov chain with stationary distributions from \citet{healy2008randomness}.
\begin{lemma}[Theorem 1.1 in \citet{healy2008randomness}]\label{lem:Markovconcentration}
Let $\bm{X}=(X_i, i\geq1)$ be a Markov chain with a stationary distribution $\phi$. Suppose that the distribution of $X_1$ is identical to the distribution of $\phi$. Then, there exists a constant $\lambda>0$ such that for any $\epsilon>0$, it holds that
\begin{equation}\label{eqn:Markovchain}
    P\left( \left|\sum_{i=1}^mX_i-\mathbb{E}\left[\sum_{i=1}^mX_i\right]\geq \sqrt{m}\cdot\epsilon\right| \right)\leq2\exp\left(-\frac{\epsilon^2(1-\lambda)}{4}\right)
\end{equation}
\end{lemma}
We also state the following lemma regarding the convexity of the pseudo-cost.
\begin{lemma}[Proposition 1 in \citet{bu2020constant}]\label{lem:convex}
We denote by $\hat{s}(\mu,Z)=s(q(\mu), z)$ where $q(\mu)=\min_{q}\{q:\mathbb{E}[s(q,Z)]\geq\mu\}$. We also denote the transformed cost function $\text{TC}(\mu)=\hat{C}_{\infty}^{\pi_{q(\mu)}}$. Suppose that the random supply function takes one of the four formulations specified in \Cref{sec:randomsupply}. Then, $\text{TC}(\mu)$ is a convex function over $[0,\bar{\mu}]$, where $\bar{\mu}$ satisfying $q(\bar{\mu})=\bar{q}$.
\end{lemma}
We finally state the following result, showing how the limiting inventory level can be bounded.
\begin{lemma}[Lundberg's Inequality]\label{lem:Lundberg}
Denote by $I_{\infty}$ as the limiting distribution of the stochastic process $I_{t+1}=(I_t+Q-D)^+$, where $Q$ and $D$ are two positive random variables. Then, there exists a constant $\rho$ such that for any $a>0$, we have
\[
P(I_{\infty}\geq a)\leq\exp(-\rho a).
\]
Moreover, $\rho$ is the adjustment coefficient of the random variable $Q-D$, which is defined as the solution to $\lambda(z)=1$,
where $\lambda(z)=\mathbb{E}[\exp(z\cdot(Q-D))]$.
\end{lemma}

\section{Missing Proofs}\label{sec:pfap}

\begin{myproof}[Proof of \Cref{lem:Lipschitz}]
From \Cref{lem:convex}, we know that the transformed cost function $\text{TC}(\mu)=\hat{C}_{\infty}^{\pi_{q(\mu)}}$ is a convex function over $\mu\in[0,\bar{\mu}]$, which is a bounded region. Thus, we know that there exists a constant $\beta'>0$ such that
\begin{equation}\label{eqn:lem201}
|\text{TC}(\mu_1)-\text{TC}(\mu_2)|\leq\beta'\cdot|\mu_1-\mu_2|,~\forall \mu_1, \mu_2\in[0,\bar{\mu}].
\end{equation}
For any $q_1, q_2\in[0,\bar{q}]$, we now denote by $\mu_1=\mathbb{E}[s(q_1,Z)]$ and $\mu_2=\mathbb{E}[s(q_2,Z)]$. Moreover, for the random supply function taking one of the four formulations specified in \Cref{sec:randomsupply}, it is direct to check that there exists a constant $\alpha'>0$ such that
\begin{equation}\label{eqn:lem202}
|\mu_1-\mu_2|\leq\alpha'\cdot|q_1-q_2|
\end{equation}
Plugging \eqref{eqn:lem202} into \eqref{eqn:lem201}, we know that
\[
|\hat{C}^{\pi_{q_1}}_{\infty}-\hat{C}^{\pi_{q_2}}_{\infty}|=|\text{TC}(\mu_1)-\text{TC}(\mu_2)|\leq\beta'\cdot|\mu_1-\mu_2|\leq\alpha'\beta'\cdot|q_1-q_2|
\]
Therefore, we prove that $\hat{C}^{\pi_q}_{\infty}$ is Lipschitz continuous over $q$ with a Lipschitz constant $\beta=\alpha'\beta'$.
\end{myproof}

\begin{myproof}[Proof of \Cref{lem:Costcoupling}]
For each epoch $n$, we denote by
\[
\mathcal{B}_n=\{ I_{\tau_{n'}}\leq \kappa_1\cdot\log T,~\tilde{I}^{a^{n*}}_{\tau_{n'}}\leq\kappa_1\cdot\log T\text{~and~}I_{\tau_{n'}+\kappa_2\cdot\max\{\log T, 2L\}}=\tilde{I}^{a^{n'*}}_{\tau_{n'}+\kappa_2\cdot\max\{\log T, 2L\}}\},~\forall n'\leq n \}.
\]
In order to prove the lemma, it is sufficient to prove that
\begin{equation}\label{pf:boundonBn}
    P(\mathcal{B}_n)\geq 1-\frac{3n}{T^2}.
\end{equation}
We prove \eqref{pf:boundonBn} by using induction on the epoch $n$.

Clearly, when $n=1$, we have that $P(I_{\tau_1}=0)=1$. From \Cref{lem:Lundberg}, there exists a constant $\kappa_1>0$ such that
\[
P(\tilde{I}^{a^{1*}}_{\tau_1}\leq\kappa_1\cdot\log T)\geq 1-\frac{1}{T^2}
\]
by noting that the distribution of $\tilde{I}^{a^{1*}}_{\tau_1}$ is identical to the distribution of $I^{a^{1*}}_{\infty}$.

Now conditioning on the event $\{\tilde{I}^{a^{1*}}_{\tau_1}\leq\kappa_1\cdot\log T\}$, we proceed to bound the probability that event $\{ I_{\tau_{1}+\kappa_2\cdot\max\{\log T, 2L\}}=\tilde{I}^{a^{1*}}_{\tau_{1}+\kappa_2\cdot\max\{\log T, 2L\}} \}$ happens. Note that the evolution of the stochastic process $I_t$ in \eqref{eqn:Itevolve} is identical to the evolution of the stochastic process $\tilde{I}^{a^{1*}}_t$ in \eqref{eqn:tIevolve} for $t=\tau_1,\dots,\tau_2-1$. Therefore, it is clear to see that the event $\{I_{\tau_1+\kappa_2\cdot\max\{\log T, 2L\}}=\tilde{I}^{a^{1*}}_{\tau_1+\kappa_2\cdot\max\{\log T, 2L\}}\}$ happens as long as
\begin{equation}\label{pf:equal010}
I_{t}=\tilde{I}^{a^{1*}}_{t}=0,~~\text{for~some~}t=\tau_1,\dots,\tau_1+\kappa_2\cdot\max\{\log T, 2L\}.
\end{equation}
We note that $I_{\tau_1}\leq\tilde{I}^{a^{1*}}_{\tau_1}$ implies that $I_{t}\leq \tilde{I}^{a^{1*}}_t$ for all $t=\tau_1,\dots,\tau_1+\kappa_2\cdot\max\{\log T, 2L\}$. From the non-negativity of $I_t$ and $\tilde{I}^{a^{1*}}_t$, we have that \eqref{pf:equal010} holds as long as $\tilde{I}_t^{a^{1*}}=0$ for some $t=\tau_1,\dots,\tau_1+\kappa_2\cdot\max\{\log T, 2L\}$. As a result, a sufficient condition for \eqref{pf:equal010} to hold is that
\[
\sum_{t=\tau_1}^{\tau_1+\kappa_2\cdot\max\{\log T,2L\}}D_t-s(a^{1*},Z_t)\geq\kappa_1\cdot\log T
\]
Note that $D_t-s(a^{1*},Z_t)$ are i.i.d. random varibles for $t=\tau_1,\dots,\tau_1+\kappa_2\cdot\max\{\log T, 2L\}$. We also denote by $\delta=\mathbb{E}[D]-\mathbb{E}[s(\bar{q},Z)]$. Following Hoeffding's inequality (\Cref{lem:Hoeffding}), we have that
\[
P\left( \sum_{t=\tau_1}^{\tau_1+\kappa_2\cdot\max\{\log T,2L\}}D_t-s(a^{1*},Z_t)\geq\kappa_1\cdot\log T \right)\geq1-\exp(-\frac{2(\delta\kappa_2\max\{\log T,2L\}-\kappa_1\log T)^2}{\kappa_2\cdot\max\{\log T, 2L\}\cdot\bar{D}})\geq1-\frac{1}{T^2}
\]
where $\kappa_2\geq\max\{\frac{2\kappa_1}{\delta},\frac{4\bar{D}}{\delta^2}\}\geq\max\{\frac{2\kappa_1\log T}{\delta\cdot\max\{\log T,2L\}}, \frac{4\bar{D}\log T}{\delta^2\cdot\max\{\log T,2L\}}\}$.

The above derivation implies that
\[\begin{aligned}
P(\mathcal{B}_1)&=P\left( \mathcal{B}_1\mid\tilde{I}^{a^{1*}}_{\tau_1}\leq\kappa_1\cdot\log T\right)\cdot P(\tilde{I}^{a^{1*}}_{\tau_1}\leq\kappa_1\cdot\log T)\\
&=P\left( \sum_{t=\tau_1}^{\tau_1+\kappa_2\cdot\max\{\log T,2L\}}D_t-s(a^{1*},Z_t)\geq\kappa_1\cdot\log T \right)\cdot P(\tilde{I}^{a^{1*}}_{\tau_1}\leq\kappa_1\cdot\log T)\\
&\geq(1-\frac{1}{T^2})\cdot(1-\frac{1}{T^2})\geq1-\frac{2}{T^2}\geq1-\frac{3}{T^2}
\end{aligned}\]
Therefore, we prove \eqref{pf:boundonBn} for $n=1$.

Now assume that \eqref{pf:boundonBn} holds for epoch $n-1$. We consider epoch $n$. Clearly, from the definition of the stochastic process $\tilde{I}_t^{a^{n*}}$ in \eqref{eqn:tIevolve}, $\tilde{I}_t^{a^{n*}}$ refreshes when $t=\tau_n$. As a result, the distribution of $\tilde{I}_{\tau_n}^{a^{n*}}$ is independent of the event $\mathcal{B}_{n-1}$ and is identical to the distribution of $I_{\infty}^{a^{n*}}$, which implies that
\begin{equation}\label{pf:new01}
P(\tilde{I}^{a^{n*}}_{\tau_n}\leq\kappa_1\cdot\log T\mid\mathcal{B}_{n-1})=P(\tilde{I}^{a^{n*}}_{\tau_n}\leq\kappa_1\cdot\log T)\geq 1-\frac{1}{T^2}
\end{equation}
where the second inequality follows from \Cref{lem:Lundberg}. Moreover, conditioning on $\mathcal{B}_{n-1}$, since $I_{\tau_n-1}$ couples with $\tilde{I}_{\tau_n-1}^{a^{(n-1)*}}$, we have that
\begin{equation}\label{pf:new02}
P(I_{\tau_n-1}\leq\kappa_1\cdot\log T\mid\mathcal{B}_{n-1})=P(\tilde{I}_{\tau_n-1}^{a^{(n-1)*}}\leq\kappa_1\cdot\log T\mid\mathcal{B}_{n-1})=P(\tilde{I}_{\tau_n-1}^{a^{(n-1)*}}\leq\kappa_1\cdot\log T)/P(\mathcal{B}_{n-1})
\end{equation}
Note that the distribution of $\tilde{I}_{\tau_n-1}^{a^{(n-1)*}}$ is identical to the distribution of $I_{\infty}^{a^{(n-1)*}}$, which implies that
\[
P(\tilde{I}_{\tau_n-1}^{a^{(n-1)*}}\leq\kappa_1\cdot\log T)=P(I_{\infty}^{a^{(n-1)*}}\leq\kappa_1\cdot\log T)\geq1-\frac{1}{T^2}
\]
Therefore, by noting that $P(\mathcal{B}_{n-1})\leq1$, from \eqref{pf:new02}, we have that
\begin{equation}\label{pf:new03}
    P(I_{\tau_n-1}\leq\kappa_1\cdot\log T\mid\mathcal{B}_{n-1})\geq1-\frac{1}{T^2}
\end{equation}
From \eqref{pf:new01}, \eqref{pf:new03} and the union bound, we have that
\begin{equation}\label{pf:new04}
    P(I_{\tau_n-1}\leq\kappa_1\cdot\log T\text{~and~}\tilde{I}^{a^{n*}}_{\tau_n}\leq\kappa_1\cdot\log T \mid\mathcal{B}_{n-1})\geq1-\frac{2}{T^2}
\end{equation}
As a result, conditioning on $\mathcal{B}_{n-1}$, we know that
\begin{equation}\label{pf:new05}
I_{\tau_n+L}\leq L\cdot\bar{D}+\kappa_1\cdot\log T\text{~and~}\tilde{I}^{a^{n*}}_{\tau_n+L}\leq L\cdot\bar{D}+\kappa_1\cdot\log T
\end{equation}
happens with a probability at least $1-\frac{2}{T^2}$.
It is clear to see that the event $\{I_{\tau_n+\max\{\log T, 2L\}}=\tilde{I}^{a^{n*}}_{\tau_n+\max\{\log T, 2L\}}\}$ happens as long as
\begin{equation}\label{pf:equal0}
I_{t}=\tilde{I}^{a^{n*}}_{t}=0
\end{equation}
for some $t=\tau_n+L,\dots,\tau_n+\max\{\log T, 2L\}$.

Suppose that $I_{\tau_n}\leq\tilde{I}^{a^{n*}}_{\tau_n}$ (resp. $I_{\tau_n}\geq\tilde{I}^{a^{n*}}_{\tau_n}$), from the evolution of the stochastic process in \eqref{eqn:Itevolve} and \eqref{eqn:tIevolve}, we have that $I_{t}\leq\tilde{I}^{a^{n*}}_{t}$ (resp. $I_{t}\geq\tilde{I}^{a^{n*}}_{t}$) for any $t=\tau_n+L,\dots,\tau_n+\max\{\log T, 2L\}$. Given that $I_t$ and $\tilde{I}^{a^{n*}}_t$ must be non-negative (from definition), we conclude that if $I_{\tau_n}\leq\tilde{I}^{a^{n*}}_{\tau_n}$ (resp. $I_{\tau_n}\geq\tilde{I}^{a^{n*}}_{\tau_n}$), then \eqref{pf:equal0} happens as long as $\tilde{I}^{a^{n*}}_{t}=0$ (resp. $I_{t}=0$). Thus, a sufficient condition for \eqref{pf:equal0} to happen is that
\[
\sum_{t=\tau_n+L}^{\tau_n+\max\{\log T,2L\}}D_t-s(a^{n*},Z_t)\geq L\cdot\bar{D}+\kappa_1\cdot\log T
\]
Since $D_t-s(a^{n*},Z_t)$ are i.i.d. random variable for $t=\tau_n+L, \dots, \tau_n+\max\{\log T,2L\}$, we denote by $\delta_n=\mathbb{E}_{D\sim F}[D]-\mathbb{E}_{Z\sim G}[s(a^{n*},Z)]\geq\delta$. Following Hoeffding's inequality (\Cref{lem:Hoeffding}), we have that
\[\begin{aligned}
P\left(\sum_{t=\tau_n+L}^{\tau_n+\kappa_2\max\{\log T,2L\}}D_t-s(a^{n*},Z_t)\geq L\cdot\bar{D}+\kappa_1\cdot\log T\right)&\geq1-\exp(-\frac{2(\kappa_2\max\{\log T,2L\}-L(\bar{D}+1)-\kappa_1\log T)^2}{\kappa_2\max\{\log T,2L\}-L})\\
&\geq 1-\frac{1}{T^2}
\end{aligned} \]
where $\kappa_2\geq\max\{4,2(\bar{D}+1+\kappa_1)\}\geq\max\{\frac{4\log T}{\max\{\log T,2L\}}, 2(\bar{D}+1+\kappa_1)\}$. Therefore, we have that
\[
P\left( I_{\tau_n+\kappa_2\max\{\log T, 2L\}}=\tilde{I}^{a^{n*}}_{\tau_n+\kappa_2\max\{\log T, 2L\}}\mid \mathcal{B}_{n-1}\text{~and~}\eqref{pf:new05}\text{~happens} \right)\geq 1-\frac{1}{T^2}.
\]
Combining \eqref{pf:new04} and the induction hypothesis that $P(\mathcal{B}_{n-1})\geq1-\frac{3(n-1)}{T^2}$, we have that
\[\begin{aligned}
P(\mathcal{B}_n)=&P(\mathcal{B}_{n-1})\cdot P(I_{\tau_n-1}\leq\kappa_1\cdot\log T\text{~and~}\tilde{I}^{a^{n*}}_{\tau_n}\leq\kappa_1\cdot\log T \mid\mathcal{B}_{n-1})\\
&\cdot P\left( I_{\tau_n+\kappa_2\max\{\log T, 2L\}}=\tilde{I}^{a^{n*}}_{\tau_n+\kappa_2\max\{\log T, 2L\}}\mid \mathcal{B}_{n-1}\text{~and~}\eqref{pf:new05}\text{~happens} \right)\\
\geq&(1-\frac{3(n-1)}{T^2})\cdot(1-\frac{2}{T^2})\cdot(1-\frac{1}{T^2})\geq(1-\frac{3(n-1)}{T^2})\cdot(1-\frac{3}{T^2})\\
\geq&1-\frac{3n}{T^2}
\end{aligned}\]
which completes our proof of the induction of \eqref{pf:boundonBn} for each epoch $n$. Therefore, our proof of the lemma is completed.
\end{myproof}

\begin{myproof}[Proof of \Cref{lem:Gapcost}]
Clearly, from \eqref{eqn:Reintercost}, it is enough to compare the value of $I_t$ and $\tilde{I}_t^{a^{n*}}$ for each epoch $n$ and each period $t$ in the epoch $n$. Note that we identify an event $\mathcal{B}$ in \Cref{lem:Costcoupling} that $I_t$ and $\tilde{I}_t^{a^{n*}}$ couple with each other. We consider two situations where $\mathcal{B}$ happens or $\mathcal{B}$ not happens.

Case 1: We now assume that $\mathcal{B}$ happens. Then, we know that for each epoch $n\in[N]$ and each $t=\tau_n+\kappa_2\cdot\max\{\log T, 2L\},\dots,\tau_{n+1}-1$, the value of $I_t$ and $\tilde{I}_t^{a^{n*}}$ are identical. Therefore, only when $t=\tau_n,\dots,\tau_n+\kappa_2\cdot\max\{\log T, 2L\}$, the value of $I_t$ and $\tilde{I}_t^{a^{n*}}$ can be different. Moreover, note that the evolution of $I_t$ in \eqref{eqn:Itevolve} is the same as the evolution of $\tilde{I}_t^{a^{n*}}$ in \eqref{eqn:tIevolve}, except that the intial value $I_{\tau_n}$ is different from $\tilde{I}_{\tau_n}^{a^{n*}}$. We know that the gap between $I_t$ and $\tilde{I}_t^{a^{n*}}$ can only become smaller. Therefore, we get that
\begin{equation}\label{pf:lem304}
    |I_t-\tilde{I}_t^{a^{n*}}|\leq|I_{\tau_n}-\tilde{I}_{\tau_n}^{a^{n*}}|\leq\kappa_1\cdot\log T
\end{equation}
where the last inequality follows from the condition in the event $\mathcal{B}$. We have that
\begin{equation}\label{pf:lem303}
\begin{aligned}
    \left|\mathbb{E}\left[\sum_n\sum_{t=\tau_n}^{\tau_{n+1}-1}(I_t-\tilde{I}_t^{a^{n*}})\mid \mathcal{B}\right]\right|&\leq\sum_n\sum_{t=\tau_n}^{\tau_n+\kappa_2\cdot\max\{\log T,2L\}}\mathbb{E}[|I_t-\tilde{I}_t^{a^{n*}}|\mid\mathcal{B}]\\
    &\leq\sum_n \kappa_1\cdot\log T\cdot \kappa_2\cdot\max\{\log T,2L\}\\
    &=N\cdot \kappa_1\kappa_2\log T\cdot\max\{\log T,2L\}
\end{aligned}
\end{equation}

Case 2: We now assume that $\mathcal{B}$ does not happen. Clearly, a direct upper bound on both $I_t$ and $\tilde{I}_t^{a^{n*}}$ is that
\[
I_t\leq\bar{D}\cdot t \text{~and~}\mathbb{E}[\tilde{I}_t^{a^{n*}}\mid\mathcal{B}^c]\leq \bar{D}\cdot t
\]
Therefore, we have that
\begin{equation}\label{pf:lem304}
\left|\mathbb{E}\left[\sum_n\sum_{t=\tau_n}^{\tau_{n+1}-1}(I_t-\tilde{I}_t^{a^{n*}})\mid \mathcal{B}^c\right]\right|\leq\bar{D}\cdot\sum_{t=1}^Tt\leq\bar{D}\cdot T^2
\end{equation}
However, from \Cref{lem:Costcoupling}, we know that $P(\mathcal{B}^c)\leq\frac{3N}{T^2}$. As a result, combining \eqref{pf:lem303} and \eqref{pf:lem304}, we get that
\[\begin{aligned}
\left|\mathbb{E}\left[\sum_n\sum_{t=\tau_n}^{\tau_{n+1}-1}(I_t-\tilde{I}_t^{a^{n*}})\right]\right|&\leq\left|\mathbb{E}\left[\sum_n\sum_{t=\tau_n}^{\tau_{n+1}-1}(I_t-\tilde{I}_t^{a^{n*}})\mid \mathcal{B}\right]\right|\cdot P(\mathcal{B})+\left|\mathbb{E}\left[\sum_n\sum_{t=\tau_n}^{\tau_{n+1}-1}(I_t-\tilde{I}_t^{a^{n*}})\mid \mathcal{B}^c\right]\right|\cdot P(\mathcal{B}^c)\\
&\leq \left|\mathbb{E}\left[\sum_n\sum_{t=\tau_n}^{\tau_{n+1}-1}(I_t-\tilde{I}_t^{a^{n*}})\mid \mathcal{B}\right]\right|+\left|\mathbb{E}\left[\sum_n\sum_{t=\tau_n}^{\tau_{n+1}-1}(I_t-\tilde{I}_t^{a^{n*}})\mid \mathcal{B}^c\right]\right|\cdot\frac{3N}{T^2}\\
&\leq N\cdot \kappa_1\kappa_2\log T\cdot\max\{\log T,2L\}+3N\bar{D}
\end{aligned}\]
which completes our proof.
\end{myproof}

\begin{myproof}[Proof of \Cref{lem:confidencecoupling}]
The proof generalizes the proof of \Cref{lem:Costcoupling}.
For each epoch $n$, we denote by
\[
\mathcal{C}_n=\{ I^a_{\tau_{n'}}\leq \kappa_1\cdot\log T,~\tilde{I}^a_{\tau_{n'}}\leq\kappa_1\cdot\log T\text{~and~}I^a_{\tau_{n'}+\kappa_2\cdot\max\{\log T, 2L\}}=\tilde{I}^{a^{n'*}}_{\tau_{n'}+\kappa_2\cdot\max\{\log T, 2L\}}\},~\forall a\in\mathcal{A}_{n'}, ~\forall n'\leq n \}.
\]
In order to prove the lemma, it is sufficient to prove that
\begin{equation}\label{pf:lem4boundonBn}
    P(\mathcal{C}_n)\geq 1-\frac{3(K+1)n}{T^2}.
\end{equation}
We prove \eqref{pf:lem4boundonBn} by using induction on the epoch $n$.

Clearly, when $n=1$, we have that $P(I^a_{\tau_1}=0)=1$ for all $a\in\mathcal{A}_1$. From \Cref{lem:Lundberg}, there exists a constant $\kappa_1>0$ such that for each $a\in\mathcal{A}_1$, it holds that
\[
P(\tilde{I}^{a}_{\tau_1}\leq\kappa_1\cdot\log T)\geq 1-\frac{1}{T^2}
\]
by noting that the distribution of $\tilde{I}^{a}_{\tau_1}$ is identical to the distribution of $I^{a}_{\infty}$ for each $a\in\mathcal{A}_1$.

Now conditioning on the event $\{\tilde{I}^{a}_{\tau_1}\leq\kappa_1\cdot\log T,~\forall a\in\mathcal{A}_1\}$, we proceed to bound the probability that event $\{ I^a_{\tau_{1}+\kappa_2\cdot\max\{\log T, 2L\}}=\tilde{I}^{a}_{\tau_{1}+\kappa_2\cdot\max\{\log T, 2L\}},~\forall a\in\mathcal{A}_1 \}$ happens. Note that the evolution of the stochastic process $I_t^a$ in \eqref{eqn:aevolve} is identical to the evolution of the stochastic process $\tilde{I}^{a}_t$ in \eqref{eqn:tildeIat} for $t=\tau_1,\dots,\tau_2-1$. Therefore, for each $a\in\mathcal{A}_1$, it is clear to see that the event $\{I^a_{\tau_1+\kappa_2\cdot\max\{\log T, 2L\}}=\tilde{I}^{a}_{\tau_1+\kappa_2\cdot\max\{\log T, 2L\}}\}$ happens as long as
\begin{equation}\label{pf:lem401}
I^a_{t}=\tilde{I}^{a}_{t}=0,~~\text{for~some~}t=\tau_1,\dots,\tau_1+\kappa_2\cdot\max\{\log T, 2L\}.
\end{equation}
We note that $I^a_{\tau_1}\leq\tilde{I}^{a}_{\tau_1}$ implies that $I^a_{t}\leq \tilde{I}^{a}_t$ for all $t=\tau_1,\dots,\tau_1+\kappa_2\cdot\max\{\log T, 2L\}$. From the non-negativity of $I_t^a$ and $\tilde{I}^{a}_t$, we have that \eqref{pf:lem401} holds as long as $\tilde{I}_t^{a}=0$ for some $t=\tau_1,\dots,\tau_1+\kappa_2\cdot\max\{\log T, 2L\}$. As a result, a sufficient condition for \eqref{pf:lem401} to hold for a $a\in\mathcal{A}_1$ is that
\[
\sum_{t=\tau_1}^{\tau_1+\kappa_2\cdot\max\{\log T,2L\}}D_t-s(a,Z_t)\geq\kappa_1\cdot\log T
\]
Note that $D_t-s(a,Z_t)$ are i.i.d. random varibles for $t=\tau_1,\dots,\tau_1+\kappa_2\cdot\max\{\log T, 2L\}$. We also denote by $\delta=\mathbb{E}[D]-\mathbb{E}[s(\bar{q},Z)]$. Following Hoeffding's inequality (\Cref{lem:Hoeffding}), we have that
\[
P\left( \sum_{t=\tau_1}^{\tau_1+\kappa_2\cdot\max\{\log T,2L\}}D_t-s(a,Z_t)\geq\kappa_1\cdot\log T \right)\geq1-\exp(-\frac{2(\delta\kappa_2\max\{\log T,2L\}-\kappa_1\log T)^2}{\kappa_2\cdot\max\{\log T, 2L\}\cdot\bar{D}})\geq1-\frac{1}{T^2}
\]
where $\kappa_2\geq\max\{\frac{2\kappa_1}{\delta}, \frac{4\bar{D}}{\delta^2}\}\geq\max\{\frac{2\kappa_1\log T}{\delta\cdot\max\{\log T,2L\}}, \frac{4\bar{D}\log T}{\delta^2\cdot\max\{\log T,2L\}}\}$.

The above derivation implies that
\[\begin{aligned}
P(\mathcal{B}_1)&=P\left( \mathcal{B}_1\mid\tilde{I}^{a}_{\tau_1}\leq\kappa_1\cdot\log T,~\forall a\in\mathcal{A}_1\right)\cdot P(\tilde{I}^{a}_{\tau_1}\leq\kappa_1\cdot\log T, ~\forall a\in\mathcal{A}_1)\\
&=P\left( \sum_{t=\tau_1}^{\tau_1+\kappa_2\cdot\max\{\log T,2L\}}D_t-s(a,Z_t)\geq\kappa_1\cdot\log T,~\forall a\in\mathcal{A}_1 \right)\cdot P(\tilde{I}^{a^{1*}}_{\tau_1}\leq\kappa_1\cdot\log T,~\forall a\in\mathcal{A}_1)\\
&\geq(1-\frac{K+1}{T^2})\cdot(1-\frac{K+1}{T^2})\geq1-\frac{2(K+1)}{T^2}\geq1-\frac{3(K+1)}{T^2}
\end{aligned}\]
where the first inequality follows from the union bound by noting that $|\mathcal{A}_1|\leq K+1$. Therefore, we prove \eqref{pf:lem4boundonBn} for $n=1$.

Now assume that \eqref{pf:lem4boundonBn} holds for epoch $n-1$. We consider epoch $n$. Clearly, from the definition of the stochastic process $\tilde{I}_t^{a}$ in \eqref{eqn:tildeIat}, $\tilde{I}_t^{a}$ refreshes when $t=\tau_n$. As a result, the distribution of $\tilde{I}_{\tau_n}^{a}$ is independent of the event $\mathcal{C}_{n-1}$ and is identical to the distribution of $I_{\infty}^{a}$, which implies that
\begin{equation}\label{pf:lem4new01}
P(\tilde{I}^{a}_{\tau_n}\leq\kappa_1\cdot\log T\mid\mathcal{C}_{n-1})=P(\tilde{I}^{a}_{\tau_n}\leq\kappa_1\cdot\log T)\geq 1-\frac{1}{T^2},~~\forall a\in\mathcal{A}_n
\end{equation}
where the second inequality follows from \Cref{lem:Lundberg}. Moreover, conditioning on $\mathcal{B}_{n-1}$, since $I_{\tau_n-1}^a$ couples with $\tilde{I}_{\tau_n-1}^{a}$ for each $a\in\mathcal{A}_n\subset\mathcal{A}_{n-1}$, we have that
\begin{equation}\label{pf:lem4new02}
P(I^a_{\tau_n-1}\leq\kappa_1\cdot\log T\mid\mathcal{C}_{n-1})=P(\tilde{I}_{\tau_n-1}^{a}\leq\kappa_1\cdot\log T\mid\mathcal{C}_{n-1})=P(\tilde{I}_{\tau_n-1}^{a}\leq\kappa_1\cdot\log T)/P(\mathcal{C}_{n-1}),~\forall a\in\mathcal{A}_n
\end{equation}
Note that the distribution of $\tilde{I}_{\tau_n-1}^{a}$ is identical to the distribution of $I_{\infty}^{a}$, which implies that
\[
P(\tilde{I}_{\tau_n-1}^{a}\leq\kappa_1\cdot\log T)=P(I_{\infty}^{a}\leq\kappa_1\cdot\log T)\geq1-\frac{1}{T^2},~\forall a\in\mathcal{A}_n
\]
Therefore, by noting that $P(\mathcal{C}_{n-1})\leq1$, from \eqref{pf:lem4new02}, we have that
\begin{equation}\label{pf:lem4new03}
    P(I_{\tau_n-1}^a\leq\kappa_1\cdot\log T\mid\mathcal{C}_{n-1})\geq1-\frac{1}{T^2},~\forall a\in\mathcal{A}_n.
\end{equation}
From \eqref{pf:lem4new01}, \eqref{pf:lem4new03} and the union bound, we have that
\begin{equation}\label{pf:lem4new04}
    P(I_{\tau_n-1}\leq\kappa_1\cdot\log T\text{~and~}\tilde{I}^{a^{n*}}_{\tau_n}\leq\kappa_1\cdot\log T,~\forall a\in\mathcal{A}_n \mid\mathcal{C}_{n-1})\geq1-\frac{2(K+1)}{T^2}
\end{equation}
where we note that $|\mathcal{A}_n|\leq K+1$.
As a result, conditioning on $\mathcal{C}_{n-1}$, we know that
\begin{equation}\label{pf:lem4new05}
I_{\tau_n+L}^a\leq L\cdot\bar{D}+\kappa_1\cdot\log T\text{~and~}\tilde{I}^{a}_{\tau_n+L}\leq L\cdot\bar{D}+\kappa_1\cdot\log T,~\forall a\in\mathcal{A}_n
\end{equation}
happens with a probability at least $1-\frac{2(K+1)}{T^2}$.
It is clear to see that for each $a\in\mathcal{A}_n$, the event $\{I^a_{\tau_n+\max\{\log T, 2L\}}=\tilde{I}^{a}_{\tau_n+\max\{\log T, 2L\}}\}$ happens as long as
\begin{equation}\label{pf:lem4equal0}
I^a_{t}=\tilde{I}^{a}_{t}=0
\end{equation}
for some $t=\tau_n+L,\dots,\tau_n+\max\{\log T, 2L\}$.

For each $a\in\mathcal{A}_n$, suppose that $I^a_{\tau_n}\leq\tilde{I}^{a}_{\tau_n}$ (resp. $I^a_{\tau_n}\geq\tilde{I}^{a}_{\tau_n}$), from the evolution of the stochastic process in \eqref{eqn:aevolve} and \eqref{eqn:tildeIat}, we have that $I_{t}^a\leq\tilde{I}^{a}_{t}$ (resp. $I^a_{t}\geq\tilde{I}^{a}_{t}$) for any $t=\tau_n+L,\dots,\tau_n+\max\{\log T, 2L\}$. Given that $I^a_t$ and $\tilde{I}^{a}_t$ must be non-negative (from definition), we conclude that if $I^a_{\tau_n}\leq\tilde{I}^{a}_{\tau_n}$ (resp. $I^a_{\tau_n}\geq\tilde{I}^{a}_{\tau_n}$), then \eqref{pf:equal0} happens as long as $\tilde{I}^{a}_{t}=0$ (resp. $I^a_{t}=0$). Thus, a sufficient condition for \eqref{pf:lem4equal0} to happen is that
\[
\sum_{t=\tau_n+L}^{\tau_n+\max\{\log T,2L\}}D_t-s(a,Z_t)\geq L\cdot\bar{D}+\kappa_1\cdot\log T
\]
Since $D_t-s(a,Z_t)$ are i.i.d. random variable for $t=\tau_n+L, \dots, \tau_n+\max\{\log T,2L\}$, we denote by $\delta_{n,a}=\mathbb{E}_{D\sim F}[D]-\mathbb{E}_{Z\sim G}[s(a,Z)]\geq\delta$. Following Hoeffding's inequality (\Cref{lem:Hoeffding}), for each $a\in\mathcal{A}_n$, we have that
\[\begin{aligned}
P\left(\sum_{t=\tau_n+L}^{\tau_n+\kappa_2\max\{\log T,2L\}}D_t-s(a,Z_t)\geq L\cdot\bar{D}+\kappa_1\cdot\log T\right)&\geq1-\exp(-\frac{2(\kappa_2\max\{\log T,2L\}-L(\bar{D}+1)-\kappa_1\log T)^2}{\kappa_2\max\{\log T,2L\}-L})\\
&\geq 1-\frac{1}{T^2}
\end{aligned} \]
where $\kappa_2\geq\max\{4,2(\bar{D}+1+\kappa_1)\}\geq\max\{\frac{4\log T}{\max\{\log T,2L\}}, 2(\bar{D}+1+\kappa_1)\}$. Therefore, from the union bound, we have that
\[
P\left( I^a_{\tau_n+\kappa_2\max\{\log T, 2L\}}=\tilde{I}^{a}_{\tau_n+\kappa_2\max\{\log T, 2L\}},~\forall a\in\mathcal{A}_n\mid \mathcal{C}_{n-1}\text{~and~}\eqref{pf:lem4new05}\text{~happens} \right)\geq 1-\frac{K+1}{T^2}.
\]
Combining \eqref{pf:lem4new04} and the induction hypothesis that $P(\mathcal{C}_{n-1})\geq1-\frac{3(K+1)(n-1)}{T^2}$, we have that
\[\begin{aligned}
P(\mathcal{C}_n)=&P(\mathcal{C}_{n-1})\cdot P(I^a_{\tau_n-1}\leq\kappa_1\cdot\log T\text{~and~}\tilde{I}^{a}_{\tau_n}\leq\kappa_1\cdot\log T,~\forall a\in\mathcal{A}_n \mid\mathcal{C}_{n-1})\\
&\cdot P\left( I^a_{\tau_n+\kappa_2\max\{\log T, 2L\}}=\tilde{I}^{a}_{\tau_n+\kappa_2\max\{\log T, 2L\}},~\forall a\in\mathcal{A}_n\mid \mathcal{C}_{n-1}\text{~and~}\eqref{pf:lem4new05}\text{~happens} \right)\\
\geq&(1-\frac{3(K+1)(n-1)}{T^2})\cdot(1-\frac{2(K+1)}{T^2})\cdot(1-\frac{K+1}{T^2})\\
\geq&(1-\frac{3(K+1)(n-1)}{T^2})\cdot(1-\frac{3(K+1)}{T^2})\\
\geq&1-\frac{3(K+1)n}{T^2}
\end{aligned}\]
which completes our proof of the induction of \eqref{pf:lem4boundonBn} for each epoch $n$. Therefore, our proof of the lemma is completed.
\end{myproof}

\begin{myproof}[Proof of \Cref{lem:highprob}]
We first show that for each epoch $n\in[N]$ and each action $a\in\mathcal{A}_n$, we can use the average value of $\tilde{I}^a_t$ for $t=\tau_n+\kappa_2\cdot\max\{\log T,2L\}$ to $\tau_{n+1}-1$ to approximate the value of $\mathbb{E}[I^{\pi_a}_{\infty}]$, where the length of the confidence interval can be given by $\gamma_n$.

Clearly, $\{\tilde{I}^a_t\}_{t=\tau_n+\kappa_2\cdot\max\{\log T,2L\}}^{\tau_{n+1}-1}$ forms a Markov chain. We denote by $\bm{I}$ a vector such that
\[
\bm{I}=(\tilde{I}_t^a, \forall t=\tau_n+\kappa_2\cdot\max\{\log T,2L\},\dots,\tau_{n+1}-1)
\]
We apply \Cref{lem:Markovconcentration} to derive a concentration bound for $\bm{I}$. To be specific, for each epoch $n\leq N-1$, we regard  $\tilde{I}^a_{\tau_n+\kappa_2\cdot\max\{\log T,2L\}+i}$ as $X_i$ for $i=1,\dots,\tau_{n+1}-1-\tau_n-\kappa_2\cdot\max\{\log T,2L\}$. Clearly, $\bm{I}$ is a Markov chain with stationary distributions and satisfies the conditions in \Cref{lem:Markovconcentration}.

We now denote $m=\tau_{n+1}-\tau_n-1-\kappa_2\cdot\max\{\log T,2L\}$.
Then, from \Cref{lem:Markovconcentration}, there exists a constant $\lambda$ such that for any $\epsilon>0$, it holds
\[
P\left( \left|\sum_{i=1}^m\tilde{I}^a_{\tau_n+\kappa_2\cdot\max\{\log T,2L\}+i}-\mathbb{E}\left[ \sum_{i=1}^m\tilde{I}^a_{\tau_n+\kappa_2\cdot\max\{\log T,2L\}+i} \right] \right|\geq\sqrt{m}\cdot\epsilon\right)\leq2\exp\left(-\frac{\epsilon^2(1-\lambda)}{4}\right).
\]
We now set $\epsilon=\frac{\gamma_n\cdot\sqrt{m}}{2}$. Then, we have
\begin{equation}\label{pf:lem601}
\exp\left(-\frac{\epsilon^2(1-\lambda)}{4}\right)\leq\exp(-\frac{(1-\lambda)m\gamma_n^2}{16})
\end{equation}
We proceed to give a lower bound on $m\gamma_n^2$, which will imply an upper bound for \eqref{pf:lem601}. Note that
\[
m=\tau_{n+1}-\tau_n-1-\kappa_2\cdot\max\{\log T,2L\}=\kappa_2\cdot(\max\{\frac{1}{\gamma_n^2}\cdot\log T,3L\}-\max\{\log T,2L\})
\]
If $\frac{1}{\gamma_2}\cdot\log T\geq 3L$, then we have
\[
\max\{\frac{1}{\gamma_n^2}\cdot\log T,3L\}-\max\{\log T,2L\}=\frac{1}{\gamma_n^2}\cdot\log T-\max\{\log T,2L\}\geq\frac{1}{3\gamma_n^2}\cdot\log T
\]
If $\frac{1}{\gamma_2}\cdot\log T< 3L$, then we have
\[
\max\{\frac{1}{\gamma_n^2}\cdot\log T,3L\}-\max\{\log T,2L\}=L\geq \frac{1}{3\gamma_n^2}\cdot\log T
\]
Therefore, it holds that
\begin{equation}\label{pf:lem602}
m=\kappa_2\cdot(\max\{\frac{1}{\gamma_n^2}\cdot\log T,3L\}-\max\{\log T,2L\})\geq\frac{\kappa_2}{3\gamma_n^2}\cdot\log T
\end{equation}
Plugging \eqref{pf:lem602} into \eqref{pf:lem601}, we have
\[
\exp\left(-\frac{\epsilon^2(1-\lambda)}{4}\right)\leq\exp\left(-\frac{(1-\lambda)\kappa_2\log T}{12}\right)\leq\frac{1}{T^2}
\]
where $\kappa_2\geq 24/(1-\lambda)$. We have
\begin{equation}\label{pf:lem603}
\begin{aligned}
    &P\left( \left|\sum_{i=1}^m\tilde{I}^a_{\tau_n+\kappa_2\cdot\max\{\log T,2L\}+i}-\mathbb{E}\left[ \sum_{i=1}^m\tilde{I}^a_{\tau_n+\kappa_2\cdot\max\{\log T,2L\}+i} \right] \right|
    \geq m\cdot\frac{\gamma_n}{2}\right)\\
    =& P\left( \left|\sum_{i=1}^m\tilde{I}^a_{\tau_n+\kappa_2\cdot\max\{\log T,2L\}+i}-m\cdot\mathbb{E}\left[ I^{\pi_a}_{\infty} \right] \right|
    \geq m\cdot\frac{\gamma_n}{2}\right)  \\
    \leq& \frac{2}{T^2}
\end{aligned}
\end{equation}
where the second inequality follows from the fact that the distribution of $\tilde{I}^a_t$ is identical to the distribution of $I^{\pi_a}_{\infty}$. Moreover, from Hoeffding's inequality (\Cref{lem:Hoeffding}), it holds that
\begin{align}
P\left(\left|\sum_{t=\tau_n+\kappa_2\cdot\max\{\log T,2L\}}^{\tau_{n+1}-1} s(a,Z_{t}))-m\cdot\mathbb{E}[s(a,Z)]\right|\geq m\cdot\frac{\gamma_n}{2} \right)&\leq2\exp(-\frac{m\gamma_n^2}{2\bar{D}^2})\leq 2\exp(-\frac{\kappa_2\log T}{6\bar{D}^2})\nonumber\\
&\leq\frac{2}{T^2}\label{pf:lem604}
\end{align}
where the second inequality follows from \eqref{pf:lem602} and the third inequality follows from $\kappa_2\geq12\bar{D}^2$. Therefore, conditional on the event $\mathcal{C}$ happens, we have that
\[
P\left(\left|\tilde{C}^a_n-\hat{C}^{\pi_a}_{\infty}\right|\leq(h+b)\cdot\frac{\gamma_n}{2}\mid\mathcal{C}\right)\geq1-\frac{4}{T^2}
\]
which implies that (from union bound over all $a\in\mathcal{A}$ and all $n\leq N-1$)
\[
P(\mathcal{E}\mid\mathcal{C})=P\left(\{ |\tilde{C}^a_n-\hat{C}^{\pi_a}_{\infty}|\leq(h+b)\cdot\frac{\gamma_n}{2},~\forall a\in \mathcal{A}_n, \forall 1\leq n\leq N-1 \}\right)\geq1-\frac{4(K+1)N}{T^2}
\]
From \Cref{lem:confidencecoupling}, we know that $P(\mathcal{C})\geq1-\frac{3(K+1)N}{T^2}$. Therefore, we have that
\[
P(\mathcal{E})=P(\mathcal{E}\mid\mathcal{C})\cdot P(\mathcal{C})\geq1-\frac{7(K+1)N}{T^2}
\]
which completes our proof.
\end{myproof}

\end{APPENDICES}

\end{document}